\documentclass[12pt, a4paper]{article}
\usepackage[latin1]{inputenc}
\usepackage{amsmath}
\usepackage[english]{babel}
\usepackage{indentfirst}
\usepackage{amstext}
\usepackage{amsfonts}
\usepackage{textcomp}
\usepackage{color}
\usepackage{amssymb}
\usepackage[dvips]{graphicx}
\usepackage{setspace}
\usepackage[all]{xy}

\newcommand {\demo}{\hskip -0.6cm{\bf Proof.  }}
\newcommand {\fim}{\hfill{$\square$}\vskip 1pc}

\newcommand {\R}{\mathbb{R}}
\newcommand {\N}{\mathbb{N}}
\newcommand {\Z}{\mathbb{Z}}

\newcommand{\supp}{\text{supp}}

\newcommand{\Per}{\mathrm{Per}}

\newcommand{\id}{\mathrm{id}}

\newtheorem{teorema}{Theorem}[section]
\newtheorem{lema}[teorema]{Lemma}

\newtheorem{definicao}[teorema]{Definition}
\newtheorem{proposicao}[teorema]{Proposition}
\newtheorem{exemplo}[teorema]{Example}
\newtheorem{rmk}[teorema]{Remark}
\newtheorem{ntn}[teorema]{Notation}

\begin{document}

\onehalfspace

\title{Branching Systems for Higher-Rank Graph C*-algebras}
\author{Daniel Gon\c{c}alves\footnote{Partially supported by CNPq, Brazil.}, Hui Li\footnote{Partially supported by Research Center for Operator Algebras of East China Normal University; partially supported by Science and Technology Commission of Shanghai Municipality (STCSM), grant No. 13dz2260400; and partially supported by an NSERC Discovery grant of Dilian Yang. The corresponding author.} \ and Danilo Royer}
\date{15 Mar 2017}
\maketitle

AMS 2010 MSC: 46L05, 37A55

Keywords: $k$-graph; $k$-graph C*-algebra; branching system; single-vertex $k$-graph C*-algebra; faithful representation.

\begin{abstract}
We define branching systems for finitely aligned higher-rank graphs. From these we construct concrete representations of higher-rank graph C*-algebras on Hilbert spaces. We prove a generalized Cuntz-Krieger uniqueness theorem for periodic single-vertex $2$-graphs. We use this result to give a sufficient condition under which representations of periodic single-vertex $2$-graph C*-algebras arising from branching systems are faithful.
\end{abstract}

\section{Introduction}

Higher-rank graphs, or $k$-graphs, are combinatorial objects that generalize directed graphs. In \cite{KumjianPask:NYJM00} Kumijian and Pask introduced higher-rank graph C*-algebras for row-finite higher-rank graphs without sources, as generalizations of graph algebras and the higher-rank Cuntz-Krieger algebras constructed by Robertson and Steger \cite{Robertson-Steger}. Since then, driven by the fact that higher-rank graph C*-algebras includes a larger class than graph C*-algebras, while still can be studied via combinatorial methods, intense research has been done in the subject, see \cite{MR3150172, MR3150171, DavidsonYang, MR2511133, MR3315580, KPSW, RaeburnSimsEtAl:PEMS03, RaeburnSimsEtAl:JFA04, MR3392275}, for example.

Branching systems arise in disciplines such as random walk, symbolic dynamics and scientific computing (see for example \cite{Dev07, HG09, SaSt96}). More recently, stimulated by Bratteli-Jorgensen's work connecting representations of the Cuntz algebra arising from iterated function systems with wavelets (see \cite{MR1465320, BJ}), a large number of papers have studied representations of graph algebras from branching systems (see \cite{chen, GLR1, GLR, GR4, MR2903145, GR3, GR, MR2848777, HR}, etc). Farsi, Gillaspy, Kang and Packer have studied connections of representations of finite higher-rank graphs C*-algebras arising from semibranching function systems with wavelets, KMS states (see \cite{MR3404559, FGKP}, etc.).

It is our intention to connect the theory of higher-rank graph C*-algebras with the branching system theory. Notice that when developing the theory of higher-rank graphs some additional hypotheses are usually assumed, such as finiteness, row finiteness, local convexity or finite alignment. Of these, finite alignment is the most general one, and we try to study the branching system theory of higher-rank graphs in this generality as much as we can. As the paper goes on, to obtain interesting results, we reduce the generality to row-finite higher-rank graphs without sources. Eventually we restrict to single-vertex $2$-graphs, which have been studied in depth by Davidson and Yang (see \cite{DavidsonYang}).

The structure of the paper is as follows. Section 2 is devoted to recalling the material on higher-rank graph C*-algebras. In Section 3 we define branching systems for finitely aligned higher-rank graphs. Using the space of boundary paths of a higher-rank graph, we build a branching system associated to any finitely aligned higher-rank graph. We then show how branching systems induce representations of higher-rank graph C*-algebras, which generalizes results in \cite{MR2903145}.
in Section 4 we look into some examples of higher-rank graphs and build branching systems on $\R$ for these graphs, which including higher-rank graphs that are not row-finite. It is usually not easy to decide if a representation of a higher-rank graph C*-algebra is faithful. Therefore in Section 5 we focus on studying periodic single-vertex $2$-graphs, and we aim to provide a sufficient condition for representations induced from branching systems to be faithful. To do so, we firstly extend the general Cuntz-Krieger uniqueness theorem proved by Jonathan Brown, Gabriel Nagy, and Sarah Reznikoff in \cite[Theorem~7.10]{MR3150172}, in the same spirit of Wojciech Szyma{\'n}ski's result for graph algebras (see \cite[Theorem~1.2]{Szyma'nski:IJM02}) and the author's result for ultragraph algebras (see \cite[Theorem~7.4]{GLR1}). We finish the section by building branching systems for periodic single-vertex $2$-graphs such that the  induced representations are faithful.


\section{Preliminaries}

Throughout this paper, the notation $\mathbb{N}$ stands for the set of all nonnegative integers; the notation $\mathbb{N}_+$ stands for the set of all positive integers; and all measure spaces are assumed to be $\sigma$-finite.


In this section, we recall the definition of $k$-graph C*-algebras from \cite{KumjianPask:NYJM00, RaeburnSimsEtAl:PEMS03, RaeburnSimsEtAl:JFA04}.






\begin{definicao}[{\cite[Definition~1.1]{KumjianPask:NYJM00}}]
Let $k \in \mathbb{N}_+$. A small category $\Lambda$ is called a \emph{$k$-graph} if there exists a functor $d:\Lambda \to \mathbb{N}^k$ satisfying the \emph{factorization property}, that is, for $\mu\in\Lambda, n,m \in \mathbb{N}^k$ with $d(\mu)=n+m$, there exists unique $\nu,\alpha \in \Lambda$ such that $d(\nu)=n,d(\alpha)=m, s(\nu)=r(\alpha), \mu=\nu\alpha$. The functor $d$ is called the degree map of $\Lambda$.

Let $(\Lambda_1,d_1), (\Lambda_2,d_2)$ be two $k$-graphs. A functor $f:\Lambda_1 \to \Lambda_2$ is called a \emph{morphism} if $d_2 \circ f=d_1$.
\end{definicao}

Throughout this paper, all $k$-graphs are assumed to be countable.

\begin{exemplo}[{\cite[Page~211]{RaeburnSimsEtAl:JFA04}}]
Let $k \in \mathbb{N}_+$ and let $n \in (\mathbb{N} \cup \{\infty\})^{k}$. Define $\Omega_{k,n}:=\{(p,q) \in \mathbb{N}^k \times \mathbb{N}^k:p \leq q \leq n \}$. For $(p,q), (q,m) \in \Omega_{k,n}$, define $(p,q) \cdot (q,m):=(p,m); r(p,q):=(p,p); s(p,q):=(q,q)$; and $d(p,q):=q-p$. Then $(\Omega_{k,n},d)$ is a $k$-graph.
\end{exemplo}

\begin{ntn}[{\cite[Page~211]{RaeburnSimsEtAl:JFA04}}]
Let $k \in \mathbb{N}_+$, let $\Lambda$ be a $k$-graph. Denote by
\[
X_\Lambda:=\bigcup_{n \in (\mathbb{N} \cup \{\infty\})^{k}}\{x:\Omega_{k,n} \to \Lambda:x \text{ is a graph morphism} \}.
\]
Fix a graph morphism $x:\Omega_{k,n} \to \Lambda$ for some $n \in (\mathbb{N} \cup \{\infty\})^{k}$. For $\mu \in \Lambda$ with $s(\mu)=x(0,0)$, denote by $\mu x:\Omega_{k,d(\mu)+n} \to \Lambda$ the unique graph morphism such that $(\mu x)(0,d(\mu))=\mu, (\mu x)(d(\mu),m)=x(0,m-d(\mu))$ for all $\mathbb{N}^k \ni m \leq n$. For $\mathbb{N}^k \ni m \leq n$, denote by $\sigma^m(x):\Omega_{k,n-m} \to \Lambda$ the unique graph morphism such that $\sigma^m(x)(0,l)=x(m,m+l)$ for all $\mathbb{N}^k \ni l \leq n-m$. Moreover, for $A \subset \Lambda, B \subset X_\Lambda$, denote by $AB:=\{\mu x: \mu \in A, x \in B,s(\mu)=x(0,0)\}$.
\end{ntn}

The following lemma might be well-known, however we could not find any reference to it.

\begin{lema}\label{sigma is a local homeo}
Let $k \in \mathbb{N}_+$, let $\Lambda$ be a $k$-graph, let $\mu \in \Lambda$, and let $B \subset X_\Lambda$. Then $\sigma^{d(\mu)}:\mu B \to s(\mu) B$ is a bijection.
\end{lema}

\demo
It is straightforward to see. Indeed $\sigma^{d(\mu)}(\mu x)=x$ for all $\mu x \in \mu B$ and the inverse map of $\sigma^{d(\mu)}$ is to attach $\mu$ to the elements of $B$.

\begin{definicao}[{\cite[Definition~2.8]{RaeburnSimsEtAl:JFA04}}]
Let $k \in \mathbb{N}_+$ and let $\Lambda$ be a $k$-graph. Define
\begin{align*}
\Lambda^{\leq \infty}&:=\bigcup_{n \in (\mathbb{N} \cup \{\infty\})^{k} }\{x:\Omega_{k,n} \to \Lambda \text{ is a graph morphism } : \exists \mathbb{N}^{k} \ni n_x \leq n, \text{ s.t. } \\&\forall m \in \mathbb{N}^{k} \text{ with } n_x \leq m \leq n, \text{ we have } m_i=n_i \implies x(0,m)\Lambda^{e_i}=\emptyset  \}.
\end{align*}
$\Lambda^{\leq \infty}$ is called the \emph{boundary path space} of $\Lambda$. The range and degree maps may be extended to boundary paths $x:\Omega_{k,n} \to \Lambda$ by setting $r(x):=x(0)$ and $d(x)=n$. 
\end{definicao}

\begin{ntn}
Let $k \in \mathbb{N}_+$. Denote by $e_1,\dots,e_k$ the standard basis of $\mathbb{N}^k$. For $i \geq 1$, denote by
\begin{align*}
e_i := \begin{cases}
    e_1 &\text{ if $i=1,k+1,2k+1,3k+1\dots$;} \\
    \cdots \\
    e_k &\text{ if $i=k,k+k,2k+k,3k+k,\dots$. }
\end{cases}\end{align*}
For $n,m \in \mathbb{N}^k$, denote by $\vert n\vert:=n_1+\dots +n_k; n \lor m:=(\max\{n_i,m_i\})_{i=1}^{k}$; and $n \land m:=(\min\{n_i,m_i\})_{i=1}^{k}$. Furthermore, for $z \in \mathbb{T}^k$, denote by $z^n:=z_1^{n_1}\dots z_k^{n_k}$.
\end{ntn}

\begin{ntn}[{\cite[Definitions~2.2, 2.4]{RaeburnSimsEtAl:JFA04}}]
Let $k \in \mathbb{N}_+$ and let $\Lambda$ be a $k$-graph. For $n \in \mathbb{N}^k$, denote by $\Lambda^n:=d^{-1}(n)$. For $A,B \subset \Lambda$, define $AB:=\{\mu\nu:\mu\in A,\nu\in B, s(\mu)=r(\nu)\}$. For $\mu,\nu \in \Lambda$, define $\Lambda^{\min}(\mu,\nu):=\{(\alpha,\beta) \in \Lambda \times \Lambda:\mu\alpha=\nu\beta,d(\mu\alpha)=d(\mu)\lor d(\nu)\}$. For $v \in \Lambda^0$, a subset $E$ of $v\Lambda$ is said to be \emph{exhaustive} for $v$ if, for any $\mu \in v\Lambda$, there exists $\nu \in E$ such that $\Lambda^{\min}(\mu,\nu)\neq\emptyset$.
\end{ntn}

\begin{definicao}[{\cite[Definitions~2.2]{RaeburnSimsEtAl:JFA04}}]
Let $k \in \mathbb{N}_+$. A $k$-graph $\Lambda$ is said to be \emph{finitely aligned} if, for any $\mu,\nu \in\Lambda$, we have that $\Lambda^{\min}(\mu,\nu)$ is a finite set.
\end{definicao}

\begin{definicao}[{\cite[Definition~2.5]{RaeburnSimsEtAl:JFA04}}]\label{define k-graph alg}
Let $k \in \mathbb{N}_+$ and let $\Lambda$ be a finitely aligned $k$-graph. A \emph{Cuntz-Krieger $\Lambda$-family} in a C*-algebra $B$ is a family of partial isometries $\{S_\mu\}_{\mu \in \Lambda}$ satisfying
\begin{enumerate}
\item $\{S_v\}_{v \in \Lambda^0}$ is a family of mutually orthogonal projections;
\item $S_{\mu\nu}=S_{\mu} S_{\nu}$ if $s(\mu)=r(\nu)$;
\item $S_{\mu}^* S_{\nu}=\sum_{(\alpha,\beta)\in \Lambda^{\min}(\mu,\nu)}S_\alpha S_\beta^*$ for all $\mu,\nu \in \Lambda$; and
\item\label{CK-condition for all path} $\prod_{\mu \in E}(S_v-S_\mu S_\mu^*)=0$ for all $v \in \Lambda^0$, for all finite exhaustive set $E \subset v\Lambda$.
\end{enumerate}
The C*-algebra generated by a universal Cuntz-Krieger $\Lambda$-family, denoted by $\{s_\mu\}_{\mu \in \Lambda}$, is called the \emph{$k$-graph C*-algebra} of $\Lambda$ and is denoted by $C^*(\Lambda)$.
\end{definicao}

\begin{rmk}
By \cite[Proposition~2.12]{RaeburnSimsEtAl:JFA04}, each $s_\mu$ is nonzero.
\end{rmk}

\begin{teorema}[{\cite[Theorem~C.1]{RaeburnSimsEtAl:JFA04}}]\label{simplification of CK-family fin ali case}
Let $k \in \mathbb{N}_+$, let $\Lambda$ be a finitely aligned $k$-graph, and let $\{S_\mu:\mu \in \big(\bigcup_{i=1}^{k}\Lambda^{e_i}\big) \cup \Lambda^0\}$ be a family of partial isometries in a C*-algebra $B$ satisfying
\begin{enumerate}
\item\label{S_v mut ort proj} $\{S_v\}_{v \in \Lambda^0}$ is a family of mutually orthogonal projections;
\item\label{S_mu S_nu=S_alpha S_beta} $S_{\mu} S_{\nu}=S_{\alpha}S_\beta$ if $\mu,\nu,\alpha,\beta \in \big(\bigcup_{i=1}^{k}\Lambda^{e_i}\big) \cup \Lambda^0, \mu\nu=\alpha\beta$;
\item\label{S_mu^*S_nu} $S_{\mu}^* S_{\nu}=\sum_{(\alpha,\beta)\in \Lambda^{\min}(\mu,\nu)}S_\alpha S_\beta^*$ for all $\mu,\nu \in \big(\bigcup_{i=1}^{k}\Lambda^{e_i}\big) \cup \Lambda^0$; and
\item\label{CK condition} $\prod_{\mu \in E}(S_v-S_\mu S_\mu^*)=0$ for all $v \in \Lambda^0$, for all finite exhaustive set $E \subset v\big(\bigcup_{i=1}^{k}\Lambda^{e_i}\big)$.
\end{enumerate}
Then there exists a unique Cuntz-Krieger $\Lambda$-family $\{T_\mu\}_{\mu \in \Lambda}$ in $B$ such that $T_\mu=S_\mu$ for all $\mu \in \big(\bigcup_{i=1}^{k}\Lambda^{e_i}\big) \cup \Lambda^0$.
\end{teorema}

\begin{definicao}[{\cite[Definition~1.4]{KumjianPask:NYJM00}}]
Let $k \in \mathbb{N}_+$, let $\Lambda$ be a $k$-graph. Then $\Lambda$ is said to be \emph{row-finite} if $\vert v\Lambda^{n}\vert<\infty$ for all $v \in \Lambda^0, n \in \mathbb{N}^k$. $\Lambda$ is said to \emph{have no sources} if $v\Lambda^{n} \neq \emptyset$ for all $v \in \Lambda^0, n \in \mathbb{N}^k$.
\end{definicao}

\begin{proposicao}[{\cite[Proposition~B.1]{RaeburnSimsEtAl:JFA04}}]\label{definition of CK family row fin case}
Let $k \in \mathbb{N}_+$ and let $\Lambda$ be a row-finite $k$-graph without sources. Then a family of partial isometries $\{S_\mu\}_{\mu \in \Lambda}$ in a C*-algebra $B$ is a Cuntz-Krieger $\Lambda$-family if and only if
\begin{enumerate}
\item\label{S_v ort proj} $\{S_v\}_{v \in \Lambda^0}$ is a family of mutually orthogonal projections;
\item\label{S_munu=S_muS_nu} $S_{\mu\nu}=S_{\mu} S_{\nu}$ if $s(\mu)=r(\nu)$;
\item\label{S_mu*S_mu=S_s(mu)} $S_{\mu}^* S_{\mu}=S_{s(\mu)}$ for all $\mu \in \Lambda$; and
\item\label{CK-condition for all path row finite} $S_v=\sum_{\mu \in v \Lambda^{n}}S_\mu S_\mu^*$ for all $v \in \Lambda^0, n \in \mathbb{N}^k$.
\end{enumerate}
\end{proposicao}

\begin{rmk}
Conditions~(\ref{S_v ort proj})--(\ref{CK-condition for all path}) of Proposition~\ref{definition of CK family row fin case} are exactly the definition of a Cuntz-Krieger family for row-finite without sources $k$-graphs, as given originally by Kumjian and Pask in \cite{KumjianPask:NYJM00}.
\end{rmk}

The following proposition is a special case of \cite[Theorem~C.1]{RaeburnSimsEtAl:JFA04}.

\begin{proposicao}\label{simplification of CK-family row fin case}
Let $k \in \mathbb{N}_+$, let $\Lambda$ be a row-finite $k$-graph without sources, and let $\{S_\mu:\mu \in \big(\bigcup_{i=1}^{k}\Lambda^{e_i}\big) \cup \Lambda^0\}$ be a family of partial isometries in a C*-algebra $B$ satisfying
\begin{enumerate}
\item\label{S_v mut ort proj row fin case} $\{S_v\}_{v \in \Lambda^0}$ is a family of mutually orthogonal projections;
\item\label{S_mu S_nu=S_alpha S_beta row fin case} $S_{\mu} S_{\nu}=S_{\alpha}S_\beta$ if $\mu,\nu,\alpha,\beta \in \big(\bigcup_{i=1}^{k}\Lambda^{e_i}\big) \cup \Lambda^0, \mu\nu=\alpha\beta$;
\item\label{S_mu^*S_nu row fin case} $S_{\mu}^* S_{\mu}=S_{s(\mu)}$ for all $\mu \in \big(\bigcup_{i=1}^{k}\Lambda^{e_i}\big) \cup \Lambda^0$; and
\item\label{CK condition row fin case} $S_v=\sum_{\mu \in v \Lambda^{e_i}}S_\mu S_\mu^*$ for all $v \in \Lambda^0, i=1,\dots,k$.
\end{enumerate}
Then there exists a unique Cuntz-Krieger $\Lambda$-family $\{T_\mu\}_{\mu \in \Lambda}$ in $B$ such that $T_\mu=S_\mu$ for all $\mu \in \big(\bigcup_{i=1}^{k}\Lambda^{e_i}\big) \cup \Lambda^0$.
\end{proposicao}

\demo
By an argument similar to the one in \cite[Lemma~3.1]{KumjianPask:NYJM00}, we have that $S_{\mu}^* S_{\nu}=\sum_{(\alpha,\beta)\in \Lambda^{\min}(\mu,\nu)}S_\alpha S_\beta^*$ for all $\mu,\nu \in \big(\bigcup_{i=1}^{k}\Lambda^{e_i}\big) \cup \Lambda^0$. By the factorization property of $\Lambda$, there exists a unique family of partial isometries $\{T_\mu\}_{\mu \in \Lambda}$ in $B$ such that $T_\mu=S_\mu$, for all $\mu \in \big(\bigcup_{i=1}^{k}\Lambda^{e_i}\big) \cup \Lambda^0$, and that Conditions~(\ref{S_v ort proj}), (\ref{S_munu=S_muS_nu}) of Proposition~\ref{definition of CK family row fin case} hold. It follows easily that Conditions~(\ref{S_mu*S_mu=S_s(mu)}), (\ref{CK-condition for all path row finite}) of Proposition~\ref{definition of CK family row fin case} hold. So Proposition~\ref{definition of CK family row fin case} implies that $\{T_\mu\}_{\mu \in \Lambda}$ is the unique Cuntz-Krieger $\Lambda$-family in $B$ such that $T_\mu=S_\mu$ for all $\mu \in \big(\bigcup_{i=1}^{k}\Lambda^{e_i}\big) \cup \Lambda^0$. \fim

\begin{ntn}\label{gauge action on C*(Lambda)}
Let $k \in \mathbb{N}_+$, let $\Lambda$ be a  row-finite $k$-graph without sources. Then there exists a \emph{gauge action} which is a strongly continuous group homomorphism $\gamma:\mathbb{T}^k \to \mathrm{Aut}(C^*(\Lambda))$ such that $\alpha_z(s_\mu)=z^{d(\mu)} s_\mu$ for all $z \in \mathbb{T}^k, \mu \in \Lambda$. The \emph{fixed point algebra} is the algebra $C^*(\Lambda)^\gamma=\overline{span}\{s_\mu s_\nu^*:d(\mu)=d(\nu)\}$. The gauge action yields a faithful expectation $\Phi$ from $C^*(\Lambda)$ onto $C^*(\Lambda)^\gamma$ such that $\Phi(s_\mu s_\nu^*)=\delta_{0,d(\mu)-d(\nu)}s_\mu s_\nu^*$ for all $\mu,\nu \in \Lambda$.
\end{ntn}

\begin{definicao}[{\cite[Definition~4.3]{KumjianPask:NYJM00}}]\label{define ape condition}
Let $k \in \mathbb{N}_+$ and let $\Lambda$ be a row-finite $k$-graph without sources. Denote by $\Lambda^\infty$ the set of \emph{infinite paths}, which consists of all graph morphisms from $\Omega_{k,(\infty,\dots,\infty)}$ to $\Lambda$. Then $\Lambda$ is said to be \emph{aperiodic} if, for any $v \in \Lambda^0$, there exists $x \in v\Lambda^\infty$ such that $\sigma^{n}(x) \neq \sigma^{m}(x)$ for all $n \neq m \in \mathbb{N}^k$.
\end{definicao}

The following theorem is the Cuntz-Krieger uniqueness theorem for row-finite higher-rank graphs without sources.

\begin{teorema}[{\cite[Theorem~4.6]{KumjianPask:NYJM00}}]\label{CK uni thm k-graph}
Let $k \in \mathbb{N}_+$, let $\Lambda$ be a row-finite $k$-graph without sources, and let $\pi:C^*(\Lambda) \to B$ be a homomorphism. Suppose that $\Lambda$ is aperiodic. Then $\pi$ is injective if and only if $\pi(s_v) \neq 0$ for all $v \in \Lambda^0$.
\end{teorema}

\begin{ntn}
Let $k \in \mathbb{N}_+$, let $\Lambda$ be a row-finite $k$-graph without sources. Define the set of \emph{periodicity} of $\Lambda$ by $\Per(\Lambda):=\{d(\mu)-d(\nu):(\mu,\nu) \text{ is }$ $\text{a cycline pair } \}$. By \cite[Theorem~4.6]{MR3392275} $\Lambda$ is aperiodic if and only if $\Per(\Lambda)=\{0\}$.
\end{ntn}

\begin{definicao}[{\cite{MR3150172}}]\label{define MASA and diagonal}
Let $k \in \mathbb{N}_+$ and let $\Lambda$ be a row-finite $k$-graph without sources. A pair $(\mu,\nu) \in \Lambda\times\Lambda$ is called a \emph{cycline pair} if $s(\mu)=s(\nu)$ and $\mu x=\nu x$ for all $x \in s(\mu)\Lambda^\infty$. The C*-subalgebra $\mathcal{M}:=C^*(\{s_\mu s_\nu^*: (\mu,\nu) \text{ is a cycline pair }\})$ is called the \emph{cycline subalgebra} of $C^*(\Lambda)$. Moreover, the C*-subalgebra $\mathcal{D}:=C^*(\{s_\mu s_\mu^* :\mu \in \Lambda\})$ is called the \emph{diagonal} of $C^*(\Lambda)$.
\end{definicao}

The following theorem is the general Cuntz-Krieger uniqueness theorem for row-finite higher-rank graphs without sources.

\begin{teorema}[{\cite[Theorem~7.10]{MR3150172}}]\label{general CK uni thm k-graph}
Let $k \in \mathbb{N}_+$, let $\Lambda$ be a row-finite $k$-graph without sources, and let $\pi:C^*(\Lambda) \to B$ be a homomorphism. Then $\pi$ is injective if and only if $\pi$ is injective on $\mathcal{M}$.
\end{teorema}

\section{Branching Systems of Higher-rank Graphs}

In this section we introduce the notion of branching systems of higher-rank graphs. The branching system definition will invoke the Radon-Nikodym derivative, and we refer the reader to \cite{Royden} for background on this material.

Notice that when studying the branching systems of higher-rank graphs we always consider those graphs satisfying certain hypotheses, like finiteness, row finiteness, local convexity or finite alignment. Of these assumptions finite alignment is the most general one, and we develop the theory of branching systems for higher-rank graphs in this generality as much as we can.

\subsection{Finitely Aligned Case}

\begin{definicao}\label{branchsystem fin aligned case}
Let $k \in \mathbb{N}_+$, let $\Lambda$ be a finitely aligned $k$-graph, let $(X,\eta)$ be a measure space and let $\{R_\mu,D_v\}_{\mu\in \bigcup_{i=1}^{k}\Lambda^{e_i},v\in \Lambda^0}$ be a family of measurable subsets of $X$. Suppose that
\begin{enumerate}
\item\label{R_e cap R_f =emptyset if e neq f} $R_\mu\cap R_\nu \stackrel{\eta-a.e.}{=}\emptyset$ if $\mu \neq \nu \in \Lambda^{e_i}$ for some $1 \leq i \leq k$;
\item\label{D_v cap D_w=emptyset} $D_v \cap D_w\stackrel{\eta-a.e.}{=} \emptyset$ if $v \neq w \in \Lambda^0$;
\item\label{R_mu subset D_r(mu)} $R_\mu\stackrel{\eta-a.e.}{\subseteq}D_{r(\mu)}$ for all $\mu\in \bigcup_{i=1}^{k}\Lambda^{e_i}$;
\item\label{Radon-Nikodym derivative} for each $\mu\in \bigcup_{i=1}^{k}\Lambda^{e_i}$, there exist two measurable maps $f_\mu:D_{s(\mu)}\rightarrow R_\mu$ and $f_\mu^{-1}:R_\mu \rightarrow D_{s(\mu)}$ such that $f_\mu\circ f_\mu^{-1}\stackrel{\eta-a.e.}{=}\id_{R_\mu}, f_\mu^{-1}\circ f_\mu\stackrel{\eta-a.e.}{=}\id_{D_{s(\mu)}}$, the pushforward measure $\eta \circ f_\mu$, of $f_\mu^{-1}$ in $D_{s(\mu)}$, is absolutely continuous with respect to $\eta$ in $D_{s(\mu)}$, and the pushforward measure $\eta \circ f_\mu^{-1}$, of $f_\mu$ in $R_\mu$, is absolutely continuous with respect to $\eta$ in $R_\mu$. Denote the Radon-Nikodym derivative $d(\eta \circ f_\mu)/d\eta$ by $\Phi_{f_\mu}$ and the Radon-Nikodym derivative $d(\eta \circ f_\mu^{-1} )/d\eta$ by $\Phi_{f_\mu^{-1}}$;
\item\label{f_mu circ f_nu=f_alpha circ f_beta} for $\mu,\nu,\alpha,\beta \in \bigcup_{i=1}^{k}\Lambda^{e_i}$ with $\mu\nu=\alpha\beta$, we have $f_\mu \circ f_\nu\stackrel{\eta-a.e.}{=}f_\alpha \circ f_\beta$;
\item\label{S_mu^*S_nu=sum S_alpha S_beta^*} for $\mu, \nu \in \bigcup_{i=1}^{k}\Lambda^{e_i}$ with $r(\mu)=r(\nu)$ and $d(\mu) \neq d(\nu)$, we have $f_\mu(D_{s(\mu)} \setminus \bigcup_{(\alpha,\beta) \in \Lambda^{\min}(\mu,\nu)}$ $R_\alpha)\cap f_\nu(D_{s(\nu)} \setminus \bigcup_{(\alpha,\beta) \in \Lambda^{\min}(\mu,\nu)} R_\beta)) \stackrel{\eta-a.e.}{=} \emptyset$; and
\item\label{exhaustive condition} for any $v \in \Lambda^0$, and for any finite exhaustive set $E \subset \bigcup_{i=1}^{k}v \Lambda^{e_i}$ for $v$, we have $\bigcup_{\mu \in E}R_\mu\stackrel{\eta-a.e.}{=}D_{v}$.
\end{enumerate}
We call $\{D_v, R_\mu, f_\mu\}_{\mu\in \bigcup_{i=1}^{k}\Lambda^{e_i},v\in \Lambda^0}$ a \emph{$\Lambda$-branching system} on $(X,\eta)$.
\end{definicao}


\begin{rmk} Informally we can think of the maps $f_\mu$ as "representing" the partial isometries $S_\mu$, so that the subsets $D_{s(\mu)}$ "represent" the initial projection of $S_\mu$ and  the subsets $R_\mu$ "represent" the final projection of $S_\mu$. With this in mind, the conditions we impose on the definition of a branching system become intuitive, except condition \ref{S_mu^*S_nu=sum S_alpha S_beta^*} which we feel deserves further explaining. We will keep a rather informal tone in this remark in order to explain the intuition behind this condition. Notice that we need to rephrase condition \ref{S_mu^*S_nu} of Theorem \ref{simplification of CK-family fin ali case} as one of our conditions on a branching system. Reading it directly we would like that  $f_\mu^{-1} f_\nu |_{R_\beta} = f_\alpha f_\beta^{-1}$, for all $\beta$ such that $(\alpha,\beta) \in \Lambda^{\min}(\mu,\nu)$ (notice that since for fixed $\mu, \nu$ there exists $j$ such that if $(\alpha,\beta) \in \Lambda^{\min}(\mu,\nu)$ then $\beta \in \Lambda^j$, we have that the $\{R_\beta: (\alpha,\beta) \in \Lambda^{\min}(\mu,\nu)\}$ forms a collection of a.e. disjoint sets). But if $(\alpha,\beta) \in \Lambda^{\min}(\mu,\nu)$ then $\mu \alpha = \nu \beta$ and hence $f_\mu f_\alpha = f_\nu f_\beta$ (by condition \ref{f_mu circ f_nu=f_alpha circ f_beta}), so that $ f_\alpha = f_\mu^{-1} f_\nu f_\beta$ and $f_\mu^{-1} f_\nu |_{R_\beta} = f_\alpha f_\beta^{-1}$ is satisfied. Notice that we also desire that on $D_{s(\nu)} \setminus \bigcup_{(\alpha,\beta) \in \Lambda^{\min}(\mu,\nu)} R_\beta$ the equality $f_\mu^{-1} f_\nu |_{R_\beta} = f_\alpha f_\beta^{-1}$ not necessarily hold. This is the reason we require condition \ref{S_mu^*S_nu=sum S_alpha S_beta^*}.


\end{rmk}

In the following theorem we build a branching system associated to any finitely aligned higher-rank graph using the space of boundary paths of a higher-rank graph.

\begin{teorema}\label{existenceofabranchingsystem}
Let $k \in \mathbb{N}_+$ and let $\Lambda$ be a finitely aligned $k$-graph. Then there exists a $\Lambda$-branching system.
\end{teorema}


\demo Let $X:=\Lambda^{\leq \infty}$, and let $\eta$ be the counting measure on $X$. For $v \in \Lambda^0$, define $D_v:=v \Lambda^{\leq \infty}$. For $\mu \in \bigcup_{i=1}^{k}\Lambda^{e_i}$, define $R_\mu:=\mu \Lambda^{\leq \infty}$. By \cite[Lemma~2.11]{RaeburnSimsEtAl:JFA04}, $D_v,R_\mu$ are nonempty. It is straightforward to see that $\{R_\mu,D_v\}_{\mu\in \bigcup_{i=1}^{k}\Lambda^{e_i},v\in \Lambda^0}$ satisfies Conditions~(\ref{R_e cap R_f =emptyset if e neq f})--(\ref{R_mu subset D_r(mu)}) of Definition~\ref{branchsystem fin aligned case}.

For $\mu\in \bigcup_{i=1}^{k}\Lambda^{e_i}$, Lemma~\ref{sigma is a local homeo} yields a bijection $\sigma^{d(\mu)}:R_\mu \to D_{s(\mu)}$. Denote by $f_\mu:=(\sigma^{d(\mu)})^{-1}$.

Fix $\mu,\nu,\alpha,\beta \in \bigcup_{i=1}^{k}\Lambda^{e_i}$ with $\mu\nu=\alpha\beta$. Then for $x \in D_{s(\nu)}=D_{s(\beta)}$, we have
\[
f_\mu \circ f_\nu(x)=f_\mu(\nu x)=\mu (\nu x)=\alpha (\beta x)=f_\alpha(\beta x)=f_\alpha \circ f_\beta(x).
\]
So Condition~(\ref{f_mu circ f_nu=f_alpha circ f_beta}) of Definition~\ref{branchsystem fin aligned case} holds.

Fix $\mu, \nu \in \bigcup_{i=1}^{k}\Lambda^{e_i}$ with $r(\mu)=r(\nu)$ and $d(\mu) \neq d(\nu)$. Suppose that there exist $x \in D_{s(\mu)} \setminus \bigcup_{(\alpha,\beta) \in \Lambda^{\min}(\mu,\nu)} R_\alpha$ and $y \in D_{s(\nu)} \setminus \bigcup_{(\alpha,\beta) \in \Lambda^{\min}(\mu,\nu)} R_\beta$ such that $f_\mu(x)=f_\nu(y)$. Let $z:=\mu x=\nu y$. Then there exists $n \geq d(\mu) \lor d(\nu)$ such that $z:\Omega_{k,n} \to \Lambda$ is a graph morphism. So $z(0,d(\mu) \lor d(\nu))=z(0,d(\mu))z(d(\mu),d(\mu) \lor d(\nu))=\mu\alpha_0$ for some $\alpha_0 \in \Lambda^{d(\nu)}$; and $z(0,d(\mu) \lor d(\nu))=z(0,d(\nu))z(d(\nu),d(\mu) \lor d(\nu))=\nu\beta_0$ for some $\beta_0 \in \Lambda^{d(\mu)}$. Hence $(\alpha_0,\beta_0) \in \Lambda^{\min}(\mu,\nu)$ and $z=\mu\cdot \alpha_0 \cdot\sigma^{d(\mu)\lor d(\nu)}(z)=\nu \cdot \beta_0 \cdot \sigma^{d(\mu)\lor d(\nu)}(z)$. By \cite[Lemma~2.10]{RaeburnSimsEtAl:JFA04}, $x \in R_{\alpha_0},y \in R_{\beta_0}$, which is a contradiction. Therefore Condition~(\ref{S_mu^*S_nu=sum S_alpha S_beta^*}) of Definition~\ref{branchsystem fin aligned case} holds.

Fix $v \in \Lambda^0$, and fix a finite exhaustive set $E \subset \bigcup_{i=1}^{k}v \Lambda^{e_i}$ for $v$. It is straightforward to see that $\bigcup_{\mu \in E}R_\mu \subset D_{v}$. We prove the reverse inclusion. Fix a graph morphism $x: \Omega_{k,n} \to \Lambda$ in $D_v$ (notice that $n \neq 0$). Suppose that $x \notin \bigcup_{\mu \in E}R_\mu$, for a contradiction. By the definition of $D_v$, there exists $\mathbb{N}^k \ni n_x \leq n$ such that whenever $n_x \leq m \leq n, m_i=n_i$, we have $x(0,m)\Lambda^{e_i}=\emptyset$. Since $E$ is exhaustive, there exists $\mu^1 \in E$ such that $\Lambda^{\min}(\mu^1,x(0,n_x)) \neq \emptyset$. Take an arbitrary $(\alpha,\beta) \in \Lambda^{\min}(\mu^1,x(0,n_x))$. Then $d(\beta)=d(\mu^1)$. So $n_x+d(\mu^1) \leq n$. Since $x \notin \bigcup_{\mu \in E}R_\mu$ and $E$ is exhaustive, there exists $\mu^2 \in E \setminus \{\mu\}$ such that $\Lambda^{\min}(\mu^2,x(0,n_x+d(\mu^1))) \neq \emptyset$. Then $n_x+d(\mu^1)+d(\mu^2) \leq n$. Inductively, we deduce that $n_x+ \sum_{\mu \in E}d(\mu) \leq n$. Then we are not able to find any path in $\mu \in E$ such that $\Lambda^{\min}(\mu,x(0,n_x+ \sum_{\mu \in E}d(\mu)))$ because $x \notin \bigcup_{\mu \in E}R_\mu$. Hence we get a contradiction and therefore Condition~(\ref{exhaustive condition}) of Definition~\ref{branchsystem fin aligned case} holds. \fim

Before we show that a branching system induces a representation of $C^*(\Lambda)$ on $L^2(X)$ we need the following lemma.

\begin{lema}\label{chain rule}
Let $k \in \mathbb{N}_+$, let $\Lambda$ be a finitely aligned $k$-graph, and let $\{R_\mu,D_v$, $f_\mu:\mu\in \bigcup_{i=1}^{k}\Lambda^{e_i},v\in \Lambda^0\}$ be a $\Lambda$-branching system on a measure space $(X,\eta)$. Fix $\mu,\nu \in \bigcup_{i=1}^{k}\Lambda^{e_i}$ with $s(\mu)=r(\nu)$. Then $\eta \circ f_\mu \circ f_\nu, \eta \circ f_\nu, \eta$ are measures on $D_{s(\nu)}$. Furthermore, we have that $\eta \circ f_\mu \circ f_\nu$ is absolutely continuous with respect to $\eta \circ f_\nu$, and $\eta \circ f_\nu$ is absolutely continuous with respect to $\eta$. Hence
\[
d(\eta \circ f_\mu \circ f_\nu)/d(\eta)=(\Phi_{f(\mu)} \circ f_\nu) \cdot \Phi_{f_\nu}.
\]
\end{lema}

\demo
It is straightforward to see that $\eta \circ f_\mu \circ f_\nu$ is absolutely continuous with respect to $\eta \circ f_\nu$, and $\eta \circ f_\nu$ is absolutely continuous with respect to $\eta$ by Condition~\ref{Radon-Nikodym derivative} of Definition~\ref{branchsystem fin aligned case}. By the chain rule we have
\[
d(\eta \circ f_\mu \circ f_\nu)/d(\eta)=d(\eta \circ f_\mu \circ f_\nu)/d(\eta \circ f_\nu) \cdot \Phi_{f_\nu}.
\]
We show that $d(\eta \circ f_\mu \circ f_\nu)/d(\eta \circ f_\nu)=\Phi_{f(\mu)} \circ f_\nu$. For any measurable set $E \subset D_{s(\nu)}$, we have
\[
\eta \circ f_\mu \circ f_\nu(E)=\int (d(\eta \circ f_\mu \circ f_\nu)/d(\eta \circ f_\nu)) \cdot \chi_{E} \, \mathrm{d}(\eta \circ f_\nu).
\]
Let $F:=f_\nu(E)$. Then
\begin{align*}
\eta \circ f_\mu \circ f_\nu(E)&=\eta \circ f_\mu(F)
\\&=\int \Phi_{f_\mu} \chi_F \, \mathrm{d}\eta
\\&=\int (\Phi_{f_\mu} \cdot \chi_{F}) \circ f_\nu \, \mathrm{d}(\eta \circ f_\nu)
\\&=\int (\Phi_{f_\mu} \circ f_\nu) \cdot \chi_E \, \mathrm{d}(\eta \circ f_\nu).
\end{align*}
So $d(\eta \circ f_\mu \circ f_\nu)/d(\eta \circ f_\nu)=\Phi_{f(\mu)} \circ f_\nu$ and we are done.
\fim

Next we show that branching systems induce representations of higher-rank graph C*-algebras, which is a generalization of \cite[Theorem~2.2]{MR2903145}.

\begin{teorema}\label{repinducedbybranchingsystems}
Let $k \in \mathbb{N}_+$, let $\Lambda$ be a finitely aligned $k$-graph, and let $\{D_v,R_\mu,f_\mu\}_{\mu\in \bigcup_{i=1}^{k}\Lambda^{e_i},v\in \Lambda^0}$ be a $\Lambda$-branching system on a measure space $(X,\eta)$. Then there exists a unique representation $\pi:C^*(\Lambda) \to B(L^2(X,\eta))$ such that $\pi(s_\mu)(\phi)=\Phi_{f_\mu^{-1}}^{1/2}( \phi \circ f_\mu^{-1})$ and $\pi(s_v)(\phi)=\chi_{D_v}\phi$, for all $\mu\in \bigcup_{i=1}^{k}\Lambda^{e_i},v\in \Lambda^0$ and $\phi \in \mathcal{L}^2(X,\eta)$.
\end{teorema}

\demo
For $\mu \in \bigcup_{i=1}^{k}\Lambda^{e_i}$, and for $\phi \in L^2(X,\eta)$, we have
\[
\int \vert \Phi_{f_\mu^{-1}}^{1/2}(\phi \circ f_\mu^{-1})\vert^2\, \mathrm{d}\eta=\int_{R_\mu} \vert \phi\circ f_\mu^{-1}\vert^2\, \mathrm{d}(\eta \circ f_\mu^{-1})=\int_{D_{s(\mu)}} \vert \phi \vert^2\, \mathrm{d}\eta<\infty.
\]
Define $S_\mu:L^2(X,\eta) \to L^2(X,\eta)$ by $S_\mu(\phi):=\Phi_{f_\mu^{-1}}^{1/2}( \phi \circ f_\mu^{-1})$. It is straightforward to see that $S_\mu \in B(L^2(X,\eta))$. For $\phi_1,\phi_2 \in L^2(X,\eta)$, we have
\begin{align*}
\langle \phi_1,S_\mu(\phi_2)\rangle&=\int \phi_1 \cdot \Phi_{f_\mu^{-1}}^{1/2}\cdot\overline{\phi_2 \circ f_\mu^{-1}}\, \mathrm{d}\eta
\\&=\int (\Phi_{f_\mu}^{-1/2}\circ f_\mu^{-1})\cdot\phi_1 \cdot \overline{\phi_2 \circ f_\mu^{-1}}\, \mathrm{d}\eta
\\&=\int \Phi_{f_\mu}^{-1/2}\cdot(\phi_1 \circ f_\mu)\cdot\overline{\phi_2}\, \mathrm{d}(\eta \circ f_\mu)
\\&=\int \Phi_{f_\mu}^{-1} \cdot (\Phi_{f_\mu}^{1/2}\cdot(\phi_1 \circ f_\mu)\cdot\overline{\phi_2})\, \mathrm{d}(\eta \circ f_\mu)
\\&=\int \Phi_{f_\mu}^{1/2}\cdot(\phi_1 \circ f_\mu)\cdot\overline{\phi_2}\, \mathrm{d}\eta.
\end{align*}
So $S_\mu^*(\phi)= \Phi_{f_\mu}^{1/2}\cdot( \phi \circ f_\mu)$ for all $\phi \in L^2(X,\eta)$.

Notice that, for $\mu \in \bigcup_{i=1}^{k}\Lambda^{e_i}$ and $\phi \in L^2(X,\eta)$, we have $S_\mu S_\mu^*(\phi)\stackrel{\eta-a.e.}{=} \chi_{R_\mu}\phi$. So $S_\mu$ is a partial isometry.

For $v \in \Lambda^0$, define $S_v:L^2(X,\eta) \to L^2(X,\eta)$ by $S_v(\phi):=\chi_{D_v}\phi$.

We will shot that the family $\{S_\mu, S_v\}_{\mu \in \bigcup_{i=1}^{k}\Lambda^{e_i},v \in \Lambda^0}$ satisfy the conditions of Theorem~\ref{simplification of CK-family fin ali case}.

Condition~(\ref{S_v mut ort proj}) of Theorem~\ref{simplification of CK-family fin ali case} follows from Condition~(\ref{D_v cap D_w=emptyset}) of Definition~\ref{branchsystem fin aligned case}. Condition~(\ref{S_mu S_nu=S_alpha S_beta}) of Theorem~\ref{simplification of CK-family fin ali case} follows from Condition~(\ref{f_mu circ f_nu=f_alpha circ f_beta}) of Definition~\ref{branchsystem fin aligned case} and Lemma~\ref{chain rule}.

Next we check Condition~(\ref{S_mu^*S_nu}) of Theorem~\ref{simplification of CK-family fin ali case}.

Fix $\mu,\nu \in \bigcup_{i=1}^{k}\Lambda^{e_i}$.

Case $1$. $\mu=\nu$. Then $\Lambda^{\min}(\mu,\nu)=\{(s(\mu),s(\mu))\}$. Since $S_\mu^* S_\nu=S_\mu^* S_\mu=S_{s(\mu)}$, Condition~(\ref{S_mu^*S_nu}) of Theorem~\ref{simplification of CK-family fin ali case} holds.

Case $2$. $\mu\neq\nu, d(\mu)=d(\nu)=e_i$ for some $1 \leq i \leq k$. Then $\Lambda^{\min}(\mu,\nu)=\emptyset$. By Condition~(\ref{R_e cap R_f =emptyset if e neq f}) of Definition~\ref{branchsystem fin aligned case}, we have $S_\mu^* S_\nu=0$. So Condition~(\ref{S_mu^*S_nu}) of Theorem~\ref{simplification of CK-family fin ali case} holds.

Case $3$. $d(\mu) \neq d(\nu)$. Then $d(\mu)=e_i, d(\nu)=e_j$, for some $1 \leq i \neq j \leq k$. For $(\alpha,\beta) \in \Lambda^{\min}(\mu,\nu)$, we have $S_\mu S_\alpha=S_\nu S_\beta$ because we just verified Condition~(\ref{S_mu S_nu=S_alpha S_beta}) of Theorem~\ref{simplification of CK-family fin ali case}. Then $S_\alpha S_\beta^*=S_\mu^* S_\mu S_\alpha S_\beta^*=S_\mu^* S_\nu S_\beta S_\beta^*$. So
\[
\sum_{(\alpha,\beta) \in \Lambda^{\min}(\mu,\nu)}S_\alpha S_\beta^*=\sum_{(\alpha,\beta) \in \Lambda^{\min}(\mu,\nu)}S_\mu^* S_\nu S_\beta S_\beta^*.
\]
We claim that $\sum_{(\alpha,\beta) \in \Lambda^{\min}(\mu,\nu)}S_\mu^* S_\nu S_\beta S_\beta^*=S_\mu^* S_\nu$. Fix $\phi \in L^2(X,\eta)$. Then
\begin{align*}
\sum_{(\alpha,\beta) \in \Lambda^{\min}(\mu,\nu)}S_\mu^* S_\nu S_\beta S_\beta^* \phi&=\sum_{(\alpha,\beta) \in \Lambda^{\min}(\mu,\nu)}S_\mu^* S_\nu (\chi_{R_\beta}\phi)
\\&=S_\mu^* S_\nu (\chi_{\bigcup_{(\alpha,\beta) \in \Lambda^{\min}(\mu,\nu)}R_\beta} \cdot \phi)
\\&\text{ (By Condition~(\ref{R_e cap R_f =emptyset if e neq f}) of Definition~\ref{branchsystem fin aligned case}) }
\\&=\Phi_{f_\mu}^{1/2} \cdot (\Phi_{f_\nu^{-1}}^{1/2} \circ f_\mu ) \cdot (\phi \circ f_\nu^{-1} \circ f_\mu) \cdot
\\&(\chi_{\bigcup_{(\alpha,\beta) \in \Lambda^{\min}(\mu,\nu)}R_\beta}\circ f_\nu^{-1} \circ f_\mu)
\\&=\Phi_{f_\mu}^{1/2} \cdot (\Phi_{f_\nu^{-1}}^{1/2} \circ f_\mu ) \cdot (\phi \circ f_\nu^{-1} \circ f_\mu)
\\&=S_\mu^* S_\nu \phi \text{ (By Condition~\ref{S_mu^*S_nu=sum S_alpha S_beta^*} of Definition~\ref{branchsystem fin aligned case}). }
\end{align*}
So we finish proving the claim, and hence Condition~(\ref{S_mu^*S_nu}) of Theorem~\ref{simplification of CK-family fin ali case} holds.

Finally we check Condition~(\ref{CK condition}) of Theorem~\ref{simplification of CK-family fin ali case}. Fix $v \in \Lambda^0$, fix a finite exhaustive set $E \subset \bigcup_{i=1}^{k}v\Lambda^{e_i}$, and fix $\phi \in L^2(X,\eta)$. It is straightforward to see that $\prod_{\mu \in E}(S_v-S_\mu S_\mu^*)(\phi)=\prod_{\mu \in E}(\chi_{D_v}-\chi_{R_\mu})\phi$. So by Condition~(\ref{exhaustive condition}) of Definition~\ref{branchsystem fin aligned case}, we have $\prod_{\mu \in E}(\chi_{D_v}-\chi_{R_\mu})\phi=0$. Hence Condition~(\ref{CK condition}) of Theorem~\ref{simplification of CK-family fin ali case} holds. Therefore by Theorem~\ref{simplification of CK-family fin ali case} there exists a unique Cuntz-Krieger $\Lambda$-family $\{T_\mu\}_{\mu \in \Lambda}$ in $B(L^2(X,\eta))$ such that $T_\mu=S_\mu$ for all $\mu \in \big(\bigcup_{i=1}^{k}\Lambda^{e_i}\big) \cup \Lambda^0$. By the universal property of $C^*(\Lambda)$ there exists a unique representation $\pi:C^*(\Lambda) \to B(L^2(X,\eta))$ such that $\pi(s_\mu)=T_\mu$ for all $\mu \in \Lambda$. \fim


\subsection{Row-finite Without Sources Case}

In this subsection, we simplify the definition of branching systems for row-finite $k$-graphs without sources.

\begin{definicao}\label{branchsystem row fin without sources case}
Let $k \in \mathbb{N}_+$, let $\Lambda$ be a row-finite $k$-graph without sources, let $(X,\eta)$ be a measure space and let $\{R_\mu,D_v\}_{\mu\in \bigcup_{i=1}^{k}\Lambda^{e_i},v\in \Lambda^0}$ be a family of measurable subsets of $X$. Suppose that
\begin{enumerate}
\item\label{R_e cap R_f =emptyset if e neq f row fin case} $R_\mu\cap R_\nu \stackrel{\eta-a.e.}{=}\emptyset$ if $\mu \neq \nu \in \Lambda^{e_i}$ for some $1 \leq i \leq k$;
\item\label{D_v cap D_w=emptyset row fin case} $D_v \cap D_w\stackrel{\eta-a.e.}{=} \emptyset$ if $v \neq w \in \Lambda^0$;
\item\label{Radon-Nikodym derivative row fin case} for each $\mu\in \bigcup_{i=1}^{k}\Lambda^{e_i}$, there exist two measurable maps $f_\mu:D_{s(\mu)}\rightarrow R_\mu$ and $f_\mu^{-1}:R_\mu \rightarrow D_{s(\mu)}$ such that $f_\mu\circ f_\mu^{-1}\stackrel{\eta-a.e.}{=}\id_{R_\mu}, f_\mu^{-1}\circ f_\mu\stackrel{\eta-a.e.}{=}\id_{D_{s(\mu)}}$, the pushforward measure $\eta \circ f_\mu$, of $f_\mu^{-1}$ in $D_{s(\mu)}$, is absolutely continuous with respect to $\eta$ in $D_{s(\mu)}$, and the pushforward measure $\eta \circ f_\mu^{-1}$, of $f_\mu$ in $R_\mu$, is absolutely continuous with respect to $\eta$ in $R_\mu$. Denote the Radon-Nikodym derivative $d(\eta \circ f_\mu)/d\eta$ by $\Phi_{f_\mu}$, and the Radon-Nikodym derivative $d(\eta \circ f_\mu^{-1} )/d\eta$ by $\Phi_{f_\mu^{-1}}$;
\item\label{f_mu circ f_nu=f_alpha circ f_beta row fin case} for $\mu,\nu,\alpha,\beta \in \bigcup_{i=1}^{k}\Lambda^{e_i}$ with $\mu\nu=\alpha\beta$, we have $f_\mu \circ f_\nu\stackrel{\eta-a.e.}{=}f_\alpha \circ f_\beta$;
\item\label{exhaustive condition row fin case} for $v \in \Lambda^0$, and for $1 \leq i \leq k$, we have $\bigcup_{\mu \in v\Lambda^{e_i}}R_\mu\stackrel{\eta-a.e.}{=}D_{v}$.
\end{enumerate}
We call $\{D_v,R_\mu,f_\mu\}_{\mu\in \bigcup_{i=1}^{k}\Lambda^{e_i},v\in \Lambda^0}$ a \emph{$\Lambda$-branching system} on $(X,\eta)$.
\end{definicao}

Next we show that, for row-finite without sources $k$-graphs, the above definition coincides with Definition~\ref{branchsystem fin aligned case}.

\begin{proposicao}\label{simplify def of branch sys row fin case}
Let $k \in \mathbb{N}_+$, let $\Lambda$ be a row-finite $k$-graph without sources, and let $(X,\eta)$ be a measure space. Suppose that $\{D_v,R_\mu,f_\mu\}$ is a $\Lambda$-branching system in the sense of Definition~\ref{branchsystem fin aligned case}. Then $\{D_v,R_\mu,f_\mu\}$ is a $\Lambda$-branching system in the sense of Definition~\ref{branchsystem row fin without sources case}. Conversely suppose that $\{D_v,R_\mu,f_\mu\}$ is a $\Lambda$-branching system in the sense of Definition~\ref{branchsystem row fin without sources case}. Then $\{D_v,R_\mu,f_\mu\}$ is a $\Lambda$-branching system in the sense of Definition~\ref{branchsystem fin aligned case}.
\end{proposicao}

\demo
Firstly suppose that $\{D_v,R_\mu,f_\mu\}$ is a $\Lambda$-branching system in the sense of Definition~\ref{branchsystem fin aligned case}. Then it is straightforward to see that Conditions~(\ref{R_e cap R_f =emptyset if e neq f row fin case})--(\ref{f_mu circ f_nu=f_alpha circ f_beta row fin case}) of Definition~\ref{branchsystem row fin without sources case} hold. For $v \in \Lambda^0$, and for $1 \leq i \leq k, v\Lambda^{e_i}$ is a finite exhaustive set for $v$. Then Condition~(\ref{exhaustive condition}) of Definition~\ref{branchsystem fin aligned case} implies Condition~(\ref{exhaustive condition row fin case}) of Definition~\ref{branchsystem row fin without sources case}. So $\{D_v,R_\mu,f_\mu\}$ is a $\Lambda$-branching system in the sense of Definition~\ref{branchsystem row fin without sources case}.

Conversely suppose that $\{D_v,R_\mu,f_\mu\}$ is a $\Lambda$-branching system in the sense of Definition~\ref{branchsystem row fin without sources case}. Then it is straightforward to see that Conditions~(\ref{R_e cap R_f =emptyset if e neq f})--(\ref{f_mu circ f_nu=f_alpha circ f_beta}) of Definition~\ref{branchsystem fin aligned case} hold. For $\mu, \nu \in \bigcup_{i=1}^{k}\Lambda^{e_i}$ with $r(\mu)=r(\nu), d(\mu) \neq d(\nu)$, for $\alpha \in s(\mu)\Lambda^{d(\nu)}$ with $\mu\alpha =\nu'\beta, \nu \neq \nu',d(\nu)=d(\nu')$, Condition~(\ref{f_mu circ f_nu=f_alpha circ f_beta row fin case}) of Definition~\ref{branchsystem row fin without sources case} implies that $f_\mu(R_\alpha)=f_\mu \circ f_\alpha(D_{s(\alpha)})=f_{\nu'} \circ f_\beta(D_{s(\beta)}) \subset R_{\nu'}$. So $f_\mu(R_\alpha) \cap f_\nu(D_{s(\nu)})=\emptyset$. So Condition~\ref{S_mu^*S_nu=sum S_alpha S_beta^*} of Definition~\ref{branchsystem fin aligned case} holds. For $v \in \Lambda^0$, for a finite exhaustive set $E \subset \bigcup_{i=1}^{k}v \Lambda^{e_i}$ for $v$, suppose that $\eta(D_v \setminus \bigcup_{\mu\in E}R_\mu)\neq 0$, for a contradiction. Since $\Lambda$ is row-finite without sources, there exists $\mu \in v \Lambda^{e_1+\cdots +e_k}$ such that $\eta((D_v \setminus \bigcup_{\mu\in E}R_\mu) \cap f_\mu(D_{s(\mu)}))\neq 0$. Since $E$ is exhaustive, there exist $\alpha \in E$ and $\beta \in \Lambda$ such that $\mu=\alpha\beta$. By Condition~(\ref{f_mu circ f_nu=f_alpha circ f_beta row fin case}) of Definition~\ref{branchsystem row fin without sources case}, we have $f_\mu(D_{s(\mu)})=f_\alpha \circ f_\beta(D_{s(\beta)}) \subset f_\alpha(D_{s(\alpha)})=R_\alpha$, which is a contradiction. Hence Condition~(\ref{exhaustive condition}) of Definition~\ref{branchsystem fin aligned case} holds. Therefore $\{D_v,R_\mu,f_\mu\}$ is a $\Lambda$-branching system in the sense of Definition~\ref{branchsystem fin aligned case}. \fim

\begin{rmk}
Proposition~\ref{simplify def of branch sys row fin case} yields that for branching systems of row-finite $k$-graphs without sources, Definition~\ref{branchsystem fin aligned case} is equivalent to Definition~\ref{branchsystem row fin without sources case}, and Definition~\ref{branchsystem row fin without sources case} has an easier formulation than Definition~\ref{branchsystem fin aligned case}. Therefore, from now on, whenever we consider branching systems of row-finite $k$-graphs without sources, we will not distinguish which definition we refer to.
\end{rmk}

\begin{ntn}
Let $k \in \mathbb{N}_+$, let $\Lambda$ be a row-finite $k$-graph without sources, let $\{D_v,R_\mu,f_\mu\}_{\mu\in \bigcup_{i=1}^{k}\Lambda^{e_i},v\in \Lambda^0}$ be a $\Lambda$-branching system on $(X,\eta)$, and let $\pi:C^*(\Lambda) \to B(L^2(X,\eta))$ be the representation obtained from Theorem~\ref{repinducedbybranchingsystems}. For $n \geq 1, \mu=\mu_1 \cdots \mu_n$, where $\mu_1,\dots, \mu_n \in \bigcup_{i=1}^{k}\Lambda^{e_i}$, denote by $f_\mu:=f_{\mu_1} \circ\cdots\circ f_{\mu_n}$ ($f_\mu$ is well-defined due to Condition~(\ref{f_mu circ f_nu=f_alpha circ f_beta row fin case}) of Definition~\ref{branchsystem row fin without sources case}); denote by $\Phi_{f_\mu}$ the Radon-Nikodym derivative $d(\eta \circ f_\mu)/d\eta$; and denote by $\Phi_{f_\mu^{-1}}$ the Radon-Nikodym derivative $d(\eta \circ f_\mu^{-1} )/d\eta$. It is straightforward to verify that $\pi(s_\mu)(\phi)=\Phi_{f_\mu^{-1}}^{1/2}( \phi \circ f_\mu^{-1}), \pi(s_\mu)^*(\phi)=\Phi_{f_\mu}^{1/2}( \phi \circ f_\mu), \pi(s_\mu^* s_\mu)(\phi)=\chi_{D_{s(\mu)}}\phi$, and $\pi(s_\mu s_\mu^*)(\phi)=\chi_{f_\mu(D_{s(\mu)})}\phi$, for all $\phi \in \mathcal{L}^2(X,\eta)$.

\end{ntn}

\subsection{Semibranching Function Systems}

Farsi, Gillaspy, Kang, and Packer in \cite{MR3404559} defined $\Lambda$-semibranching function systems for a finite $k$-graph without sources $\Lambda$ (being finite means that $\vert \Lambda^{n}\vert <\infty$ for all $n \in \mathbb{N}^k$). In this subsection we find connections between $\Lambda$-semibranching function systems and $\Lambda$ branching systems.

The following definition is inspired by \cite[Definitions~3.1, 3.2]{MR3404559}.

\begin{definicao}\label{define Lambda semibranching function system}
Let $\Lambda$ be a finite $k$-graph without sources, let $(X,\eta)$ be a measure space, let $\{\mathcal{D}_\mu,\mathcal{R}_\mu\}_{\mu \in \big(\bigcup_{i=1}^{k}\Lambda^{e_i}\big) \cup \Lambda^0}$ be a family of measurable subsets of $X$. Suppose that
\begin{enumerate}
\item for each $\mu\in \big(\bigcup_{i=1}^{k}\Lambda^{e_i}\big)\cup \Lambda^0$, there exist two measurable maps $\tau_\mu:\mathcal{D}_\mu \rightarrow \mathcal{R}_\mu$ and $\tau_\mu^{-1}:\mathcal{R}_\mu \rightarrow \mathcal{D}_\mu$ such that $\tau_\mu\circ \tau_\mu^{-1}\stackrel{\eta-a.e.}{=}\id_{\mathcal{R}_\mu}, \tau_\mu^{-1}\circ \tau_\mu\stackrel{\eta-a.e.}{=}\id_{\mathcal{D}_\mu}$, the pushforward measure $\eta \circ \tau_\mu$, of $\tau_\mu^{-1}$ in $\mathcal{D}_\mu$, is absolutely continuous with respect to $\eta$ in $\mathcal{D}_\mu$, and the pushforward measure $\eta \circ \tau_\mu^{-1}$, of $\tau_\mu$ in $\mathcal{R}_\mu$, is absolutely continuous with respect to $\eta$ in $\mathcal{R}_\mu$;
\item for $n \in \{0,e_1,\dots,e_k\}, X\stackrel{\eta-a.e.}{=}\bigcup_{\mu \in \Lambda^{n}}\mathcal{R}_\mu$;
\item for $n \in \{0,e_1,\dots,e_k\}$, for $\mu\neq\nu \in \Lambda^{n}, \mathcal{R}_\mu \cap \mathcal{R}_\nu \stackrel{\eta-a.e.}{=}\emptyset$;
\item for $v \in \Lambda^0, \tau_v\stackrel{\eta-a.e.}{=}\id, \eta(\mathcal{D}_v) >0$;
\item for $\mu \in \bigcup_{i=1}^{k}\Lambda^{e_i}$, we have $\mathcal{R}_\mu \stackrel{\eta-a.e.}{\subset} \mathcal{D}_{r(\mu)}, \mathcal{D}_{\mu}=\mathcal{D}_{s(\mu)}$;
\item for $n \in \{0,e_1,\dots,e_k\}$, define a measurable map $\tau^{n}:X \to X$ by $\tau^{n} \vert_{\mathcal{R}_\mu}:=\tau_\mu^{-1}$ for all $\mu \in \Lambda^{n}$. Then $\tau^{n} \circ \tau^{m}=\tau^{m} \circ \tau^{n}$ for all $n,m \in \{0,e_1,\dots,e_k\}$.
\end{enumerate}
We call $\{\mathcal{R}_\mu,\mathcal{D}_\mu,\tau_\mu,\tau^{n}: \mu\in \big(\bigcup_{i=1}^{k}\Lambda^{e_i}\big)\cup \Lambda^0,n \in \{0,e_1,\dots,e_k\}\}$ a \emph{partial $\Lambda$-semibranching function system} on $(X,\eta)$.
\end{definicao}

\begin{rmk}

For $\mu=\mu_1\cdots\mu_n \in \Lambda\setminus\Lambda^0$ where $\mu_1,\dots,\mu_n \in \bigcup_{i=1}^{k}\Lambda^{e_i}$, define $\mathcal{D}_\mu:=\mathcal{D}_{s(\mu)}$, define $\tau_\mu:=\tau_{\mu_1}\circ\cdots\circ \tau_{\mu_n}$, and define $\mathcal{R}_\mu:=\tau_\mu(\mathcal{D}_\mu)$. For $n=(n_1,\dots,n_k) \in \mathbb{N}^k \setminus \{0\}$, define $\tau^n:=n_1 \tau^{e_1} \circ \cdots \circ n_k \tau^{e_k}$. Then $\{\mathcal{R}_\mu,\mathcal{D}_\mu,\tau_\mu,\tau^{n}: \mu\in \Lambda,n \in  \mathbb{N}^k\}$ is in fact a \emph{$\Lambda$-semibranching function system} on $(X,\eta)$ as introduced by Farsi, Gillaspy, Kang, and Packer in \cite{MR3404559}.
\end{rmk}




\begin{proposicao}\label{encoding process}
Let $\Lambda$ be a finite $k$-graph without sources, let $(X,\eta)$ be a measure space, and let $\{D_v,R_\mu,f_\mu\}_{\mu\in \bigcup_{i=1}^{k}\Lambda^{e_i},v\in \Lambda^0}$ be a $\Lambda$-branching system on $(X,\eta)$. Suppose that $\eta(D_v)>0$ for all $v \in \Lambda^0$, and that $X=\bigcup_{v \in \Lambda^0}D_v$. For $v \in \Lambda^0$, define $\mathcal{D}_v=\mathcal{R}_v:=D_v$, define $\tau_v:\mathcal{D}_v \to \mathcal{R}_v$ to be the identity map. For $\mu \in \bigcup_{i=1}^{k}\Lambda^{e_i}$, define $\mathcal{D}_\mu:=D_{s(\mu)}$, define $\mathcal{R}_\mu:=R_\mu$, and define $\tau_\mu:=f_\mu$. For $n \in \{0,e_1,\dots,e_k\}$, define a measurable map $\tau^{n}:X \to X$ by $\tau^{n} \vert_{\mathcal{R}_\mu}:=\tau_\mu^{-1}$ for all $\mu \in \Lambda^{n}$. Then $\{\mathcal{R}_\mu,\mathcal{D}_\mu,\tau_\mu,\tau^{n}: \mu\in \big(\bigcup_{i=1}^{k}\Lambda^{e_i}\big)\cup \Lambda^0,n \in \{0,e_1,\dots,e_k\}\}$ is a partial $\Lambda$-semibranching function system on $(X,\eta)$.
\end{proposicao}

\demo
It is straightforward to see. \fim

\begin{proposicao}\label{decoding process}
Let $\Lambda$ be a finite $k$-graph without sources, let $(X,\eta)$ be a measure space, and let $\{\mathcal{R}_\mu,\mathcal{D}_\mu,\tau_\mu,\tau^{n}: \mu\in \big(\bigcup_{i=1}^{k}\Lambda^{e_i}\big)\cup \Lambda^0,n \in \{0,e_1,\dots$, $e_k\}\}$ be a partial $\Lambda$-semibranching function system on $(X,\eta)$. For $v \in \Lambda^0$, define $D_v:=\mathcal{D}_v$. For $\mu \in \bigcup_{i=1}^{k}\Lambda^{e_i}$, define $R_\mu:=\mathcal{R}_\mu$, and define $f_\mu:=\tau_\mu$. Then $\{D_v,R_\mu,f_\mu\}_{\mu\in \bigcup_{i=1}^{k}\Lambda^{e_i},v\in \Lambda^0}$ is a $\Lambda$-branching system on $(X,\eta)$.
\end{proposicao}

\demo
It is straightforward to see. \fim

\section{Examples of $\Lambda$-branching systems on $\R$ with the Lebesgue measure}\label{examplessection}

In this section we will present many examples of branching systems on $\R$. As we mentioned before, due to the large combinatorial possibilities permitted by the factorization property on a colored graph, we are not able to provide a general construction of branching systems on $\R$. Instead, in the examples, we provide an algorithmic way to build branching systems on $\R$, covering many examples of k-graphs in the literature.

\begin{exemplo}\label{ex1}
Let $\Gamma$ be the following 2-colored graph, where $\Gamma^0=\{v\}$, $\Gamma^{e_1}=\{f_1,f_2:n\in \N\}$ and $\Gamma^{e_2}=\{e\}$. This is an example in section 4 of \cite{FGKP}.

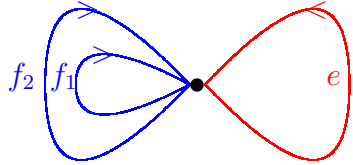
\centerline{
\setlength{\unitlength}{2cm}
\begin{picture}(1,0.7)
\color{blue}
\put(0.5,0){\qbezier(-0.05,0)(-1,-1)(-1,0)}
\put(0.5,0){\qbezier(-0.05,0)(-1,1)(-1,0)}
\put(-0.3,0.44){$>$}
\put(-0.75,0){$f_2$}
\put(0.5,0){\qbezier(-0.05,0)(-0.8,-0.4)(-0.8,0)}
\put(0.5,0){\qbezier(-0.05,0)(-0.8,0.4)(-0.8,0)}
\put(-0.2,0.142){$>$}
\put(-0.47,0){$f_1$}
\color{red}
\put(0.5,0){\qbezier(0.05,0)(1,1)(1,0)}
\put(0.5,0){\qbezier(0.05,0)(1,-1)(1,0)}
\put(1.2,0.445){$<$}
\put(1.35,0){$e$}
\color{black}
\put(0.5,0){\circle*{0.08}}
\end{picture}}
\vspace{2 cm}

There are two 2-graphs $\Lambda_2$ and $\Lambda_3$ associated to the 2-colored graph $\Gamma_2$. The factorization rules for $\Lambda_2$ are given by $$f_1 e = e f_1 \text{ and } f_2e = ef_2,$$ and the factorization rules for $\Lambda_3$ are given by $$f_1 e = e f_2 \text{ and } f_2e = ef_1.$$ For each of these 2-graphs we build a branching system below.
\end{exemplo}

We will define a branching system on $[0,1]$. Let $D_v=[0,1]$. Notice that $\{ e \}$ is exhaustive and so we must have $R_e = D_v$ (this can also be seen from Condition~\ref{exhaustive condition row fin case} of Definition~\ref{branchsystem row fin without sources case}).

To obtain a branching system for $\Lambda_2$ it is enough to define $f_{e}$ as the identity and $f_{f_1}$ and  $f_{f_2}$ as any differentiable map from $[0,1]$ to $R_{f_1}=[0,\frac{1}{2}]$, $R_{f_2}=[\frac{1}{2},1]$, respectively.

For $\Lambda_3$ we need to be a bit more careful. More precisely, the first equation of the factorization tells us that $f_e$ must take $R_{f_2}$ to $R_{f_1}$ and the second equation tells us that $R_{f_1}$ must be taken to $R_{f_2}$. Keeping the same sets as before, we define $f_e(x) = x+\frac{1}{2}$ for $x\in [0,\frac{1}{2}]$ and $f_e(x) = x-\frac{1}{2}$ for $x\in [\frac{1}{2},1]$. Furthermore, we define $f_{f_1}:[0,1]\rightarrow [0,\frac{1}{2}]$ by $f_{f_1}(x)=\frac{1}{2}x$ and the factorization now implies that $f_{f_2} = f_e f_{f_1} f_e^{-1}$ and $f_{f_2} = f_e^{-1} f_{f_1} f_e$. Since $f_e^2$ is the identity this two last equations are not contradictory and we can use it to define $f_{f_2}$ (Namely $f_{f_2}(x)= \frac{1}{2}x+ \frac{3}{4}$ if $x\in[0,\frac{1}{2}]$ and $f_{f_2}(x)= \frac{1}{2}x+ \frac{1}{4}$ if $x\in[\frac{1}{2},1]$).
\fim

\begin{exemplo}
Consider $\Lambda$ as the following 2-colored graph, where $\Lambda^0=\{v\}$, $\Lambda^{e_1}=\{g_n:n\in \N\}$ and $\Lambda^{e_2}=\{e\}$.

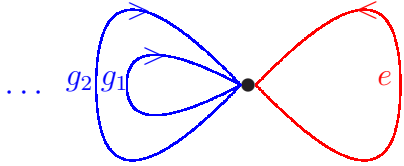
\centerline{
\setlength{\unitlength}{2cm}
\begin{picture}(1,0.7)
\put(0.5,0){\circle*{0.08}}
\color{blue}
\put(0.5,0){\qbezier(-0.05,0)(-1,-1)(-1,0)}
\put(0.5,0){\qbezier(-0.05,0)(-1,1)(-1,0)}
\put(-0.3,0.44){$>$}
\put(-0.7,0){$g_2$}
\put(0.5,0){\qbezier(-0.05,0)(-0.8,-0.4)(-0.8,0)}
\put(0.5,0){\qbezier(-0.05,0)(-0.8,0.4)(-0.8,0)}
\put(-0.2,0.142){$>$}
\put(-0.47,0){$g_1$}
\put(-1.1,-0.1){$\cdots$}
\color{red}
\put(0.5,0){\qbezier(0.05,0)(1,1)(1,0)}
\put(0.5,0){\qbezier(0.05,0)(1,-1)(1,0)}
\put(1.2,0.445){$<$}
\put(1.35,0){$e$}
\end{picture}}
\vspace{2 cm}

There are uncountably many possible factorizations in the above graph, each giving a different 2-graph. For each of these 2-graphs we build a branching system below.
\end{exemplo}

Fix a factorization and let $d:\Lambda\rightarrow \N$ be the degree map. Notice that $d(g_ie)=e_1+e_2=e_2+e_1$ and, by the factorization property, there is a unique $g_j$ such that $g_ie=eg_j$. So, we get a map $h:\N\rightarrow \N$ such that $g_ie=eg_{h(i)}$. Moreover note that $h$ is injective, because if we suppose that $h(i)=h(j)$ then $g_ie=eg_{h(i)}=eg_{h(j)}=g_je$ and then, again by the factorization property, we get $g_i=g_j$. The map $h$ is also surjective, since if $j\in \N$ then, by the factorization property, there exist a unique $i$ such that $g_ie=eg_j$.

Our goal is to define a $\Lambda$-branching system in then real interval $(0,1]$ with the Lebesgue measure. Define $D_v=(0,1]=R_e$ and $R_{g_i}=(\frac{1}{i+1},\frac{1}{i}]$, for each $i\in \N$. Now we need to define the bijective maps $\{f_{g_i}\}_{i\in \N}$ and $f_e$.

First of all, for each $i\in \N$, define the set $B_i:=\{h^n(i):n\in \Z\}$. There are two possible configurations for the $B_i$: if $h^n(i)=h^m(i)$, for some $n,m\in \Z$, then $B_i=\{i,h(i),h^2(i)...h^k(i)\}$, where the elements $h^s(i)$ are pairwise distinct and $h^{k+1}(i)=i$; if $h^n(i)\neq h^m(i)$ for each $m,n$ then $B_i=\{\cdots,h^{-2}(i), h^{-1}(i),i,h(i), h^2(i),\cdots\}$.

It is not hard to see that for $i\neq j$,  $B_i=B_j$ or $B_i\cap B_j=\emptyset$. So, by choosing an appropriated set $F\subseteq \N$ we get that $\N$ is the disjoint union $\N=\sqcup_{i\in F}B_i$.

Now we define the bijective map $f_e:D_v \rightarrow R_e$. First we define this map on each set $R_{g_i}$. Fix an $i\in F$. Suppose first that $B_i=\{i,h(i),...,h^k(i)\}$, and $h^{k+1}(i)=i$. Define, for each $n\in \{1,...,k+1\}$, the increasing linear maps $f_e:R_{g_{h^{n}(i)}}\rightarrow R_{g_{h^{n-1}(i)}}$, and piece then together to obtain the map $f_e:\bigcup\limits_{n=0}^k R_{g_{h^n(i)}}\rightarrow \bigcup\limits_{n=0}^k R_{g_{h^n(i)}}$. It follows by definition that $f_e^{k+1}$ is the identity map. For $B_i=\{...,h^{-2}(i), h^{-1}(i),i,h(i),h^2(i),\cdots\}$ define $f_e:R_{g_{h^n(i)}}\rightarrow R_{g_{h^{n-1}(i)}}$ as being the increasing linear diffeomorphism, for each $n\in \Z$. So, we get a bijective measurable map $f_e:D_v\rightarrow R_e$, with the property that, for each $i\in \N$, $f_e(R_{g_{h(i)}})=R_{g_i}$ and, moreover, if $h^{k+1}(i)=i$, for some $k\in \N$, then $f_e:\bigcup\limits_{n=0}^kR_{g_{h^n(i)}}\rightarrow \bigcup\limits_{n=0}^kR_{g_{h^n(i)}}$ is such that $f_e^{k+1}$ is the identity map (restricted to this set).

It remains to define the maps $f_{g_i}:D_v\rightarrow R_{g_i}$ for each $i\in \N$.

Let $i\in F$. If $B_i=\{i,h(i), h^2(i),...,h^k(i)\}$ with $h^{k+1}(i)=i$, define $f_{g_i}:D_v\rightarrow R_{g_i}$ as being the increasing linear diffeomorphism, and define inductively $f_{g_{h^n(i)}}=f_e^{-1}\circ f_{g_{h^{n-1}(i)}}\circ f_e$ for $n\in \{1,...,k\}$. Notice that then the equality $f_ef_{g_{h^n(i)}}=f_{g_{h^{n-1}(i)}}f_e$ holds for each $n\in \{1,...,k\}$. To see that the equality $f_ef_{g_{h^{k+1}(i)}}=f_{g_{h^k(i)}}f_e$ also holds, note that $$f_{g_{h^k(i)}}f_e=f_e^{-1}f_{g_{h^{k-1}(i)}}f_e^2=f_e^{-2}f_{g_{h^{k-2}(i)}}f_e^3=...=f_e^{-k}f_{g_{i}}f_e^{k+1}=$$
$$=f_{e}^{-k-1}f_ef_{g_i}f_e^{k+1}=f_ef_{g_i}=f_ef_{g_{h^{k+1}(i)}},$$ since $f_e^{k+1}$ is the identity map and $h^{k+1}(i)=i$.

If $B_i=\{h^{n}(i):n\in \Z\}$, with $h^{n}(i)\neq h^{m}(i)$ for each $n,m$, let $f_{g_i}:D_v\rightarrow R_{g_i}$ be the increasing linear diffeomorphism and define inductively $f_{g_{h^n(i)}}=f_e^{-1}\circ f_{g_{h^{n-1}(i)}}\circ f_e$ and $f_{g_{h^{-n}(i)}}=f_e\circ f_{g_{h^{-n+1}(i)}}\circ f_e^{-1}$ for $n\geq 1$.

It is easy to see that Conditions~(\ref{R_e cap R_f =emptyset if e neq f})--(\ref{S_mu^*S_nu=sum S_alpha S_beta^*}) of Definition \ref{branchsystem fin aligned case} are satisfied. Condition~(\ref{exhaustive condition}) also holds, because each exhaustive finite set $E\subseteq v\Lambda^{e_1}\cup v\Lambda^{e_2}$ must contain $e$ and $R_e=D_v$.
\fim

\begin{rmk} To simplify notation, and when no confusion arises, from now on we will denote the map $f_e$ associated to an edge $e$ just by $e$.
\end{rmk}

\begin{exemplo} We next turn our attention to the 2-graphs given in page 102 of \cite{RaeburnSimsEtAl:PEMS03}. Recall the 2-colored  graph, where $f$ and $h$ have degree $(0,1)$ and $k, e$ and $g$ have degree $(1,0)$:

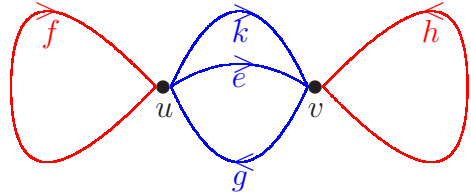
\centerline{
\setlength{\unitlength}{2cm}
\begin{picture}(1,0.7)
\put(0,0){\circle*{0.08}}
\put(-0.05,-0.2){$u$}
\put(1,0){\circle*{0.08}}
\put(0.95,-0.2){$v$}
\color{blue}
\put(0,0){\qbezier(0.05,0)(0.5,-1)(0.95,0)}
\put(0.45,-0.55){$<$}
\put(0.45,-0.65){$g$}
\put(0,0){\qbezier(0.05,0)(0.5,1)(0.95,0)}
\put(0.45,0.435){$>$}
\put(0.45,0.3){$k$}
\put(0,0){\qbezier(0.05,0)(0.5,0.3)(0.95,0)}
\put(0.45,0.1){$>$}
\put(0.45,0){$e$}
\color{red}
\put(0,0){\qbezier(-0.05,0)(-1,1)(-1,0)}
\put(0,0){\qbezier(-0.05,0)(-1,-1)(-1,0)}
\put(-0.85,0.45){$>$}
\put(-0.8,0.3){$f$}
\put(1,0){\qbezier(0.05,0)(1,1)(1,0)}
\put(1,0){\qbezier(0.05,0)(1,-1)(1,0)}
\put(1.7,0.45){$<$}
\put(1.7,0.3){$h$}
\end{picture}}
\vspace{2 cm}

There are two possible factorizations. One is $kf=hk$, $ef=he$ and $gh=fg$. For this 2-graph a branching system is obtained similarly to what we did in \ref{ex1}, defining the maps associated to the loops in the graph as the identity. We will focus in the second possible factorization, that is $$he=kf, \ hk=ef \text{ and } gh=fg.$$

\end{exemplo}

Let $D_u=[0,1]$ and $D_v=[1,2]$. Notice that the sets $\{g\}, \{f\}, \{h\}, \{e,k\}$ are exhaustive, hence $R_g=R_f=[0,1]$ and $R_h=[1,2]=R_e\cup R_k$. Let $R_e=[1, \frac{3}{2}]$ and $R_k=[\frac{3}{2},2]$. From the factorization we obtain the information on how to break up the definition of the function $h$. Let $h|_{[1, \frac{3}{2}]} \rightarrow [\frac{3}{2},2]$ be defined by $h(x)=x+\frac{1}{2}$ and $h|_{[\frac{3}{2},2]} \rightarrow [1,\frac{3}{2}]$ be defined by $h(x)=x-\frac{1}{2}$ (notice that $h^2= id$). Let $g:[1,2] \rightarrow [0,1]$ be defined by $g(x)=x-1$, $e:[0,1] \rightarrow [1,\frac{3}{2}]$ be defined by $e(x)=\frac{1}{2}x+1$ and define the remaining functions via the factorization, that is $f:=ghg^{-1}$ and $k:=h^{-1} e f=h e f = h e f^{-1}$ (notice that $f^2=id$).
\fim

\begin{exemplo} Let $\Lambda$ be the 2-graph whose 1-skeleton is given below, and where the edges $g_i$ have degree $(0,1)$ and the edges $e_i$ have degree $(1,0)$.

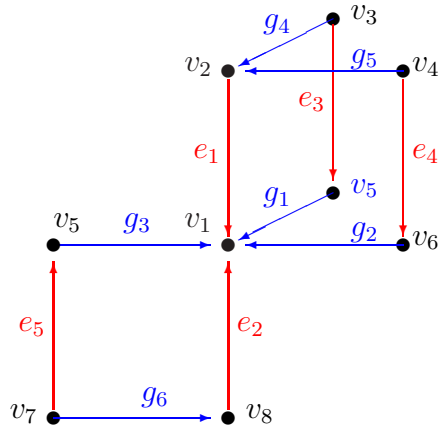
\centerline{
\setlength{\unitlength}{2.3cm}
\begin{picture}(0,1.5)
\put(0,0){\circle*{0.08}}
\put(-0.25,0.1){$v_1$}
\put(0,1){\circle*{0.08}}
\put(-0.25,1){$v_2$}
\color{red}
\put(-0.2,0.5){$e_1$}
\put(0,0.95){\vector(0,-1){0.9}}
\color{black}
\put(0.6,1.3){\circle*{0.08}}
\color{blue}
\put(0.6,1.3){\vector(-2,-1){0.54}}
\put(0.2,1.25){$g_4$}
\color{black}
\put(0.6,0.3){\circle*{0.08}}
\color{blue}
\put(0.6,0.3){\vector(-2,-1){0.54}}
\put(0.2,0.25){$g_1$}
\color{red}
\put(0.4,0.8){$e_3$}
\put(0.6,1.27){\vector(0,-1){0.9}}
\color{black}
\put(1,1){\circle*{0.08}}
\put(1,0){\circle*{0.08}}
\put(0.7,1.3){$v_3$}
\color{blue}
\put(1,1){\vector(-1,0){0.9}}
\put(1,0){\vector(-1,0){0.9}}
\put(0.7,1.05){$g_5$}
\put(0.7,0.05){$g_2$}
\put(0.7,0.3){$v_5$}
\color{red}
\put(1.05,0.5){$e_4$}
\put(1,0.95){\vector(0,-1){0.9}}
\color{black}
\put(1.05,1){$v_4$}
\put(1.05,0){$v_6$}
\put(-1,0.1){$v_5$}
\color{blue}
\put(-1,0){\vector(1,0){0.9}}
\put(-0.6,0.1){$g_3$}
\color{black}
\put(-1,0){\circle*{0.08}}
\put(-1,-1){\circle*{0.08}}
\put(0,-1){\circle*{0.08}}
\color{red}
\put(-1,-0.95){\vector(0,1){0.85}}
\put(-1.2,-0.5){$e_5$}
\color{blue}
\put(-1,-1){\vector(1,0){0.9}}
\put(-0.5,-0.9){$g_6$}
\color{black}
\put(-1.25,-1){$v_7$}
\put(0.1,-1){$v_8$}
\color{red}
\put(0,-0.95){\vector(0,1){0.85}}
\put(0.05,-0,5){$e_2$}
\end{picture}}
\vspace{2.5 cm}

\end{exemplo}

Notice that the sets $\{e_1,e_2\}$, $\{g_1,g_2,g_3\}$, $\{e_1,g_3\}$ and $\{e_2,g_1, g_2\}$ are exhaustive for $v_1$ and $\{g_4,g_5\}$ is exhaustive for $v_2$. So, take $D_{v_i}=[i-1,i]$, $R_{e_1}=[0,\frac{1}{2}]$, $R_{e_2}=[\frac{1}{2},1]$, $R_{g_3}=[\frac{1}{2},1]$, $R_{g_2}=[\frac{1}{4},\frac{1}{2}]$, $R_{g_1}=[0,\frac{1}{4}]$, $R_{g_4}=[1,\frac{3}{2}]$ and $R_{g_5}=[\frac{3}{2},2]$. Let $f_{e_1}|_{R_{g_4}}$ be the affine map onto $R_{g_1}$ and $f_{e_1}|_{R_{g_5}}$ be the affine map onto $R_{g_2}$. Define $f_{g_4}$ and $f_{e_3}$ as affine maps and let $f_{g_1}:= f_{e_1} f_{g_4} f_{g_1}^{-1}$. Analogously one define the reminder maps and sets.
\fim

\begin{exemplo} Next we construct a branching system for the 2-graph $\Lambda_2$ given on Example A.2 of \cite{RaeburnSimsEtAl:JFA04}. We reproduce a picture of the 1-skeleton below
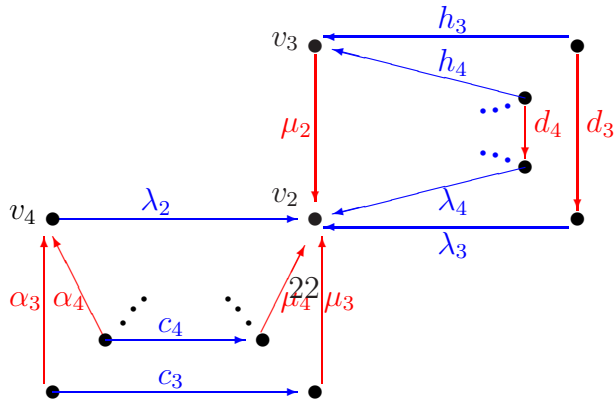
\centerline{
\setlength{\unitlength}{2.3cm}
\begin{picture}(0,1.5)
\put(0,0){\circle*{0.08}}
\put(-0.25,0.1){$v_2$}
\put(0,1){\circle*{0.08}}
\put(-0.25,1){$v_3$}
\color{red}
\put(1.2,0.65){\vector(0,-1){0.3}}
\put(-0.2,0.5){$\mu_2$}
\put(0,0.95){\vector(0,-1){0.85}}
\color{black}
\put(1.5,1){\circle*{0.08}}
\put(1.2,0.7){\circle*{0.08}}
\color{blue}
\put(1.2,0.7){\vector(-4,1){1.1}}
\color{black}
\put(1.2,0.3){\circle*{0.08}}
\color{blue}
\put(1.2,0.3){\vector(-4,-1){1.1}}
\color{black}
\put(1.5,0){\circle*{0.08}}
\color{blue}
\put(1.45,1.05){\vector(-1,0){1.4}}
\put(1.45,-0.05){\vector(-1,0){1.4}}
\put(0.7,1.1){$h_3$}
\put(0.7,0.85){$h_4$}
\put(0.7,-0.2){$\lambda_3$}
\put(0.7,0.05){$\lambda_4$}
\put(1.1,0.34){\circle*{0.03}}
\put(1.03,0.37){\circle*{0.03}}
\put(0.96,0.38){\circle*{0.03}}
\put(1.1,0.65){\circle*{0.03}}
\put(1.03,0.63){\circle*{0.03}}
\put(0.96,0.62){\circle*{0.03}}
\color{red}
\put(1.55,0.5){$d_3$}
\put(1.25,0.5){$d_4$}
\put(1.5,0.95){\vector(0,-1){0.9}}
\put(-1.2,-0.7){\vector(-1,2){0.3}}
\color{blue}
\put(-1.5,0){\vector(1,0){1.4}}
\color{black}
\put(-1.5,0){\circle*{0.08}}
\put(-1.5,-1){\circle*{0.08}}
\put(0,-1){\circle*{0.08}}
\put(-1.2,-0.7){\circle*{0.08}}
\put(-0.3,-0.7){\circle*{0.08}}
\put(-0.38,-0.6){\circle*{0.03}}
\put(-0.44,-0.52){\circle*{0.03}}
\put(-0.5,-0.46){\circle*{0.03}}
\put(-1.1,-0.6){\circle*{0.03}}
\put(-1.04,-0.52){\circle*{0.03}}
\put(-0.98,-0.46){\circle*{0.03}}
\color{red}
\put(-1.55,-0.95){\vector(0,1){0.85}}
\color{blue}
\put(-1.5,-1){\vector(1,0){1.4}}
\put(-1.2,-0.7){\vector(1,0){0.8}}
\color{black}
\put(-1.75,0){$v_4$}
\color{blue}
\put(-1,0.05){$\lambda_2$}
\put(-0.9,-0.95){$c_3$}
\put(-0.9,-0.65){$c_4$}
\color{red}
\put(-1.75,-0.5){$\alpha_3$}
\put(-1.5,-0.5){$\alpha_4$}
\put(0.05,-0.5){$\mu_3$}
\put(-0.2,-0.5){$\mu_4$}
\color{red}
\put(0.04,-0.95){\vector(0,1){0.85}}
\put(-0.31,-0.65){\vector(1,2){0.25}}
\end{picture}}
\vspace{3 cm}

In this example the edges $h_i, c_i$ and $\lambda_i$ have degree $(1,0)$ and edges $\alpha_i, \mu_i$ and $d_i$ have degree $(0,1)$.
\end{exemplo}

To construct a branching system, first we need to enumerate the edges and vertices. Respecting the labels already given in the example, we enumerate the red edges in $v_2 \Lambda_2$ by $\mu_2, \mu_3, \ldots$, the blue edges in $v_2 \Lambda_2$ by $\lambda_2, \lambda_3, \ldots$, the blue edges in $s(\mu_2)$ by $h_3, h_4, \ldots$, the red edges in $s(\lambda_2)$ by $\alpha_3, \alpha_4, \ldots$, for $i\neq 2$ we denote the red edge whose range is $s(\lambda_i)$ by $d_i$ and the blue edge whose range is $s(\mu_i)$ by $c_i$.

With the above labels, the factorization reads $$\mu_2 h_n = \lambda_n d_n \ \text{ and } \ \lambda_2 \alpha_i = \mu_i c_i.$$

As usual, to obtain a branching system, we define the sets associated to vertices as intervals. In particular, let $D_{v_2}=[0,1]$ and $D_{v_3}=[1,2]$, where $v_3= s(\mu_2)$. We will focus on defining the branching system for these vertices (once this is done it is clear how to extend it to the remaining vertices).

Notice that $\{\lambda_2, \mu_2\}$ is exhaustive. So, let $R_{\mu_2} = [0,\frac{1}{2}]$ and $R_{\lambda_2} = [\frac{1}{2},1]$. Also, let $R_{h_{i+2}}=[1+\frac{1}{2^i},1+\frac{1}{2^{i-1}}]$, $i=1,2,\ldots$ and let $R_{\lambda_{3}}=[0,\frac{1}{4}]$, $R_{\lambda_{4}}=[\frac{1}{4},\frac{1}{4}+\frac{1}{8}]$ and so on. Define $\mu_2|_{R_{h_n}}$ as the affine map onto $R_{\lambda_n}$. Also, let $h_n$ and $d_n$ be affine bijective maps. Following the factorization define, for $n\neq 2$, $\lambda_n:= \mu_2 h_n d_n^{-1}$.

Proceeding analogously one defines $\mu_i (i\neq 2), c_i, \alpha_i, \lambda_2$ and so a branching system is obtained.
\fim

\begin{exemplo} Next we construct a branching system for the 2-graph $\Lambda_3$ given on Example A.3 of \cite{RaeburnSimsEtAl:JFA04}. Differently from \cite{RaeburnSimsEtAl:JFA04} we will keep all the notations of our previous example (A.2). This example is the same as the example before, with the addition of two edges, called $\beta_3$ (of degree $(1,0)$) and $\alpha_3$ (of degree $(0,1)$) in \cite{RaeburnSimsEtAl:JFA04} (we will call $\alpha_3$ of $\nu_3$, since in our setting $\alpha_3$ is already taken). Notice that this is a particularly interesting example, since there is no finite exhaustive subset of $v_2 \Lambda$ whose range projections are orthogonal.

We reproduce a picture of the 1-skeleton below

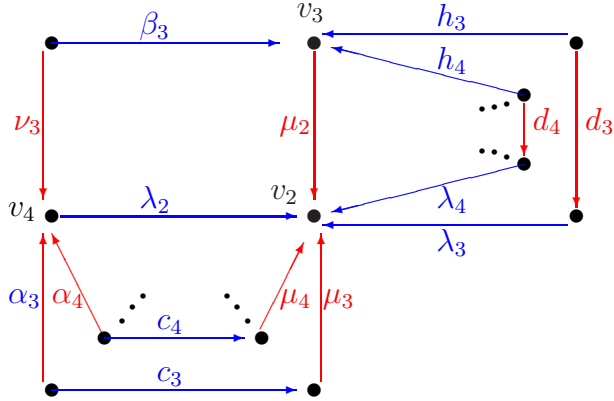
\centerline{
\setlength{\unitlength}{2.3cm}
\begin{picture}(0,1.5)
\put(0,0){\circle*{0.08}}
\put(-0.25,0.1){$v_2$}
\put(0,1){\circle*{0.08}}
\put(-0.1,1.15){$v_3$}
\color{red}
\put(1.2,0.65){\vector(0,-1){0.3}}
\put(-0.2,0.5){$\mu_2$}
\put(0,0.95){\vector(0,-1){0.85}}
\color{black}
\put(1.5,1){\circle*{0.08}}
\put(1.2,0.7){\circle*{0.08}}
\color{blue}
\put(1.2,0.7){\vector(-4,1){1.1}}
\color{black}
\put(1.2,0.3){\circle*{0.08}}
\color{blue}
\put(1.2,0.3){\vector(-4,-1){1.1}}
\color{black}
\put(1.5,0){\circle*{0.08}}
\color{blue}
\put(1.45,1.05){\vector(-1,0){1.4}}
\put(1.45,-0.05){\vector(-1,0){1.4}}
\put(0.7,1.1){$h_3$}
\put(0.7,0.85){$h_4$}
\put(0.7,-0.2){$\lambda_3$}
\put(0.7,0.05){$\lambda_4$}
\color{black}
\put(1.1,0.34){\circle*{0.03}}
\put(1.03,0.37){\circle*{0.03}}
\put(0.96,0.38){\circle*{0.03}}
\put(1.1,0.65){\circle*{0.03}}
\put(1.03,0.63){\circle*{0.03}}
\put(0.96,0.62){\circle*{0.03}}
\color{red}
\put(1.55,0.5){$d_3$}
\put(1.25,0.5){$d_4$}
\put(1.5,0.95){\vector(0,-1){0.9}}
\put(-1.2,-0.7){\vector(-1,2){0.3}}
\color{blue}
\put(-1.45,0){\vector(1,0){1.35}}
\color{black}
\put(-1.5,0){\circle*{0.08}}
\put(-1.5,-1){\circle*{0.08}}
\put(0,-1){\circle*{0.08}}
\put(-1.2,-0.7){\circle*{0.08}}
\put(-0.3,-0.7){\circle*{0.08}}
\put(-0.38,-0.6){\circle*{0.03}}
\put(-0.44,-0.52){\circle*{0.03}}
\put(-0.5,-0.46){\circle*{0.03}}
\put(-1.1,-0.6){\circle*{0.03}}
\put(-1.04,-0.52){\circle*{0.03}}
\put(-0.98,-0.46){\circle*{0.03}}
\color{red}
\put(-1.55,-0.95){\vector(0,1){0.85}}
\color{blue}
\put(-1.5,-1){\vector(1,0){1.4}}
\put(-1.2,-0.7){\vector(1,0){0.8}}
\color{black}
\put(-1.75,0){$v_4$}
\color{blue}
\put(-1,0.05){$\lambda_2$}
\put(-0.9,-0.95){$c_3$}
\put(-0.9,-0.65){$c_4$}
\put(-1.75,-0.5){$\alpha_3$}
\color{red}
\put(-1.5,-0.5){$\alpha_4$}
\put(0.05,-0.5){$\mu_3$}
\put(-0.2,-0.5){$\mu_4$}
\color{black}
\put(-1.5,1){\circle*{0.08}}
\color{blue}
\put(-1.5,1){\vector(1,0){1.3}}
\put(-1,1.05){$\beta_3$}
\color{red}
\put(-1.72,0.5){$\nu_3$}
\put(-1.55,0.95){\vector(0,-1){0.85}}
\put(0.04,-0.95){\vector(0,1){0.85}}
\put(-0.31,-0.65){\vector(1,2){0.25}}
\end{picture}}
\vspace{3 cm}

The factorization is the same as before, with one more factorization property: $$\mu_2 h_n = \lambda_n d_n \ \text{, } \ \lambda_2 \alpha_i = \mu_i c_i \ \text{ and } \ \mu_2 \beta_3 = \lambda_2 \nu_3 $$
\end{exemplo}

As before we will describe the branching system mainly at $v_2$. Notice that the only exhaustive set in $v_2$ is $\{\lambda_2, \mu_2\}$. Furthermore, the new factorization property implies that $R_{\mu_2}$ and $R_{\lambda_2}$ can not be disjoint. Let $D_{v_2}=[0,1]$, $D_{v_3}=[1,2]$, where $v_3= s(\mu_2)$, and $D_{v_4}=[2,3]$, where $v_4= s(\lambda_2)$. Define $$R_{\mu_2}=[0,\frac{3}{4}] \text{ and } R_{\lambda_2}=[\frac{1}{2},1].$$ Break $D_{v_3}$ in infinitely many intervals of positive length, namely, $\{R_{\beta_3}, R_{h_{i+2}}: i = 1, 2, \ldots\}$ and break $D_{v_4}$ in infinitely many intervals of positive length, namely, $\{R_{\nu_3}, R_{\alpha_{i+2}}: i = 1, 2, \ldots\}$. Now, break $[\frac{3}{4},1]$ in infinitely many intervals of positive length, say $R_{\mu_3}, R_{\mu_4}, \ldots$ and break $[0,\frac{1}{2}]$ in infinitely many intervals of positive length, say $R_{\lambda_3}, R_{\lambda_4}, \ldots$.

Proceeding similarly to the previous example, we define $\mu_2|_{R_{h_n}}$ as the affine map onto $R_{\lambda_n}$, $n=3,4,..$ and $\mu_2|_{R_{\beta_3}}$ as the affine map onto $R_{\mu_2}\cap R_{\lambda_2}$. We also let $h_n$ and $d_n$ be affine bijective maps and following the factorization define, for $n\neq 2$, $\lambda_n:= \mu_2 h_n d_n^{-1}$, $\lambda_2|_{R_{\nu_3}}:= \mu_2 \beta_3\nu_3^{-1}$, where $\beta_3$ and $\nu_3$ are affine bijective. The remainder of the definition of a branching systems follows analogously to above and the previous example.
\fim

\section{Faithful Representations of Periodic Single-Vertex $2$-Graphs C*-algebras via Branching Systems}

In this section we exclusively study the branching systems of periodic single-vertex $2$-graphs and we intend to find a sufficient condition for representations of periodic single-vertex $2$-graph C*-algebras, induced from branching systems, to be faithful.





First of all we recall the work of Davidson and Yang on the periodicity of single-vertex $2$-graphs in \cite{DavidsonYang} (they also studied the structure of single-vertex $k$-graph C*-algebras in \cite{MR2511133}).

\begin{teorema}[{\cite[Theorems~3.1, 3.4]{DavidsonYang}}]\label{periodicity of single-vertex 2 graphs}
Let $\Lambda$ be a single-vertex $2$-graph. Suppose that $\vert \Lambda^{e_1}\vert, \vert \Lambda^{e_2}\vert \geq 2$. Then the following are equivalent.
\begin{enumerate}
\item $\Lambda$ is periodic;
\item $\Per(\Lambda)=\mathbb{Z}(a,-b)$ for some positive integers $a,b$;
\item there exist positive integers $p,q$ with $\vert \Lambda^{e_1}\vert^p=\vert \Lambda^{e_2}\vert^q$ and a bijection $h:\prod_{i=1}^{p}\Lambda^{e_1} \to \prod_{i=1}^{q}\Lambda^{e_2}$ such that for $\mu \in \prod_{i=1}^{p}\Lambda^{e_1}, \nu \in \prod_{i=1}^{q}\Lambda^{e_2}$, we have $\mu\nu=h(\mu)h^{-1}(\nu)$ (we can identify $\prod_{i=1}^{p}\Lambda^{e_1}, \prod_{i=1}^{q}\Lambda^{e_2}$ with elements in $\Lambda$).
\end{enumerate}
\end{teorema}

\begin{ntn}\label{bijection of per 2-graph}
Let $\Lambda$ be a periodic single-vertex $2$-graph with $\vert \Lambda^{e_1}\vert, \vert \Lambda^{e_2}\vert \geq 2$. Let $(a,-b)$ be the generator of $\Per(\Lambda)$, and let $h:\prod_{i=1}^{a}\Lambda^{e_1} \to \prod_{i=1}^{b}\Lambda^{e_2}$ obtained from the above theorem. Then for each $\mu \in \prod_{i=1}^{a}\Lambda^{e_1},(\mu,h(\mu))$ is a cycline pair. By \cite[Lemma~5.3]{DavidsonYang} there is a distinguished unitary $W:=\sum_{\mu \in \prod_{i=1}^{a}\Lambda^{e_1}}s_{h(\mu)}s_{\mu}^*$ in $C^*(\Lambda)$.
\end{ntn}


\begin{lema}\label{W has full spectrum for cycline pair}
Let $\Lambda$ be a periodic single-vertex $2$-graph with $\vert \Lambda^{e_1}\vert, \vert \Lambda^{e_2}\vert\geq 2$. We inherit the notation from Notation~\ref{bijection of per 2-graph}. Then the spectrum of $W$ contains the unit circle.
\end{lema}



\demo
It is sufficient to show that $C^*(W) \cong C(\mathbb{T})$, via a unital isomorphism that identifies $W$ with the identity function on $\mathbb{T}$. By \cite[Proposition~3.11]{Katsura:CJM03} it is sufficient to show that there exists an expectation $\Phi:C^*(W) \to \mathbb{C} \cdot 1_{C^*(\Lambda)}$ such that $\Phi(W^n)=0$ for all $n \in \mathbb{Z} \setminus \{0\}$, and that $\Phi(1_{C^*(\Lambda)})=1_{C^*(\Lambda)}$. Let $\gamma$ be the gauge action on $C^*(\Lambda)$. Then $\gamma$ induces a strongly continuous homomorphism from $\mathbb{T}^2$ to $\mathrm{Aut}(C^*(W))$. Denote by $\iota:\mathbb{T} \to \mathbb{T}^2$ the embedding $z \mapsto (1,z)$. So we obtain a strongly continuous homomorphism $\gamma \circ \iota:\mathbb{T} \to \mathrm{Aut}(C^*(W))$. Therefore $\gamma\circ \iota$ yields the desired expectation $\Phi:C^*(W) \to \mathbb{C} \cdot 1_{C^*(\Lambda)}$ and hence we are done. \fim

The following theorem is an extension of the general Cuntz-Krieger uniqueness theorem of Brown-Nagy-Reznikoff (see Theorem~\ref{general CK uni thm k-graph}).

\begin{teorema}\label{general Cuntz-Krieger single vertex k-graph}
Let $\Lambda$ be a periodic single-vertex $2$-graph with $\vert \Lambda^{e_1}\vert, \vert \Lambda^{e_2}\vert\geq 2$, let $\mathcal{A}$ be a C*-algebra and let $\varphi:C^*(\Lambda) \to \mathcal{A}$ be a homomorphism. We inherit the notation from Notation~\ref{bijection of per 2-graph}. Then $\varphi$ is injective if and only if
\begin{enumerate}
\item\label{nonzero family} $\varphi(1_{C^*(\Lambda)}) \neq 0$;
\item\label{full spectrum} the spectrum of $\varphi(W)$ contains the unit circle.
\end{enumerate}
\end{teorema}

\demo
First of all suppose that $\varphi$ is injective. It is straightforward to see that Condition~(\ref{nonzero family}) holds. By Lemma~\ref{W has full spectrum for cycline pair} Condition~(\ref{full spectrum}) holds.

Conversely suppose that Conditions~(\ref{nonzero family}), (\ref{full spectrum}) hold. By Theorem~\ref{general CK uni thm k-graph}, it is sufficient to prove that $\varphi$ is injective on $\mathcal{M}$.

The faithful expectation $\Phi$ from Notation~\ref{gauge action on C*(Lambda)} restricts to a faithful expectation from $\mathcal{M}$ onto $\mathcal{D}$ satisfying that for $d \in \mathcal{D},n \in \mathbb{Z}$, if $n=0$ then $\Phi(dW^n)=d$; and if $n \neq 0$ then $\Phi(dW^n)=0$.

Since $\varphi(1_{C^*(\Lambda)}) \neq 0$, \cite[Theorem~3.1]{RaeburnSimsEtAl:JFA04} gives that $\varphi$ is injective on $C^*(\Lambda)^\gamma$. Since $\mathcal{D} \subset C^*(\Lambda)^\gamma, \varphi$ is injective on $\mathcal{D}$.

By Condition~\ref{full spectrum}, there exists an expectation $\Psi:\varphi(C^*(W)) \to \varphi(\mathbb{C} \cdot 1_{C^*(\Lambda)})$ such that for $n \in \mathbb{Z}$, if $n=0$ then $\Psi(\varphi(W^n))=\varphi(1_{C^*(\Lambda)})$; and if $n \neq 0$ then $\Psi(\varphi(W^n))=0$. As shown in the proof of \cite[Theorem~6.2]{MR3392275}, $\mathcal{M}=\overline{\mathrm{span}}\{s_\mu s_\mu^* W^n:\mu \in \Lambda,n \in \mathbb{Z}\}$ and $\mathcal{M}$ is unital abelian. We aim to construct a linear map $\Psi: \mathrm{span}\{\varphi(s_\mu s_\mu^* W^n):\mu \in \Lambda,n \in \mathbb{Z}\} \to \varphi(\mathcal{D})$ such that if $n=0$ then $\Psi(\varphi(s_\mu s_\mu^* W^n))=\varphi(s_\mu s_\mu^*)$; and if $n \neq 0$ then $\Psi(\varphi(s_\mu s_\mu^* W^n))=0$. In order to prove that $\Psi$ is well-defined, we show that it is contractive. Fix distinct $\mu_1,\dots,\mu_L \in \Lambda$ with $d(\mu_1)=\cdots=d(\mu_L)$, for $1 \leq i \leq L$, fix $\{z_{ij}\}_{j \in \mathbb{Z}} \subset \mathbb{C}$ with at most finitely many nonzero. Then we compute that
\begin{align*}
\Big\Vert \sum_{i=1}^{L}z_{i0}\varphi(s_{\mu_i} s_{\mu_i}^*) \Big\Vert&=\max_{1 \leq i \leq L}\Big\Vert z_{i0}\varphi(s_{\mu_i} s_{\mu_i}^*)  \Big\Vert
\\&\leq \max_{1 \leq i \leq L} \Big\Vert \varphi(s_{\mu_i} s_{\mu_i}^*)\sum_{j \in \mathbb{Z}}z_{ij}\varphi(W^{j}) \Big\Vert
\\&(\text{since $\Psi$ is an expectation on $\varphi(C^*(W))$})
\\&= \Big\Vert \sum_{i=1}^{L}\varphi(s_{\mu_i} s_{\mu_i}^*)\sum_{j \in \mathbb{Z}}z_{ij}\varphi(W^{j}) \Big\Vert.
\end{align*}
By Condition~(\ref{CK-condition for all path row finite}) of Proposition~\ref{definition of CK family row fin case}, every element in $\mathrm{span}\{\varphi(s_\mu s_\mu^* W^n):\mu \in \Lambda,n \in \mathbb{Z}\}$ has the form $\sum_{i=1}^{L}\varphi(s_{\mu_i} s_{\mu_i}^*)\sum_{j \in \mathbb{Z}}z_{ij}\varphi(W^{j})$. Hence we obtain a linear idempotent map $\Psi$ of norm $1$ from $\varphi(\mathcal{M})$ onto $\varphi(\mathcal{D})$. By the Tomiyama's Theorem (cf. \cite[II.6.10.2]{Blackadar:Operatoralgebras06}), $\Psi$ is an expectation. Finally by \cite[Proposition~3.11]{Katsura:CJM03}, $\varphi$ is injective on $\mathcal{M}$. So we are done. \fim






Now we present a sufficient condition for representations of periodic single-vertex $2$-graphs induced from branching systems to be faithful.

\begin{teorema}\label{a criterion of faithful rep}
Let $\Lambda$ be a periodic single-vertex $2$-graph with $\vert \Lambda^{e_1}\vert, \vert \Lambda^{e_2}\vert\geq 2$. We inherit the notation from Notation~\ref{bijection of per 2-graph}. Let $\{D_v,R_\mu,f_\mu\}_{\mu\in \bigcup_{i=1}^{k}\Lambda^{e_i}}$ be a $\Lambda$-branching system on a measure space $(X,\eta)$ such that $\eta(D_v) \neq 0$, and let $\pi:C^*(\Lambda) \to B(L^2(X,\eta))$ be the representation induced from the branching system. Suppose that for any finite subset $\mathcal{F}$ of $\mathbb{Z} \setminus \{0\}$, there exist $\mu \in \prod_{i=1}^{a}\Lambda^{e_1}$ and a measurable subset $E$ of $f_{\mu}(D_v)$ such that $\eta(E) \neq 0$ and $(f_{\mu} \circ f_{h(\mu)}^{-1})^{n}(E) \cap E\stackrel{\eta-a.e.}{=}\emptyset$ for all $n \in \mathcal{F}$. Then $\pi$ is faithful.
\end{teorema}

\demo
Since $\eta(D_v) \neq 0$, we have that $\pi(1_{C^*(\Lambda)}) \neq 0$. By \cite[Proposition~4.1]{MR3150172}, we have $\pi(s_\mu s_\mu^*)=\pi(s_{h(\mu)} s_{h(\mu)}^*)$ for all $\mu \in \prod_{i=1}^{a}\Lambda^{e_1}$. So $f_\mu(D_v)\stackrel{\eta-a.e.}{=}f_{h(\mu)}(D_v)$ for all $\mu \in \prod_{i=1}^{a}\Lambda^{e_1}; f_{\mu}(D_v) \cap f_\nu(D_v)\stackrel{\eta-a.e.}{=}\emptyset$ for distinct $\mu,\nu \in \prod_{i=1}^{a}\Lambda^{e_1}$; and $D_v\stackrel{\eta-a.e.}{=}\bigcup_{\mu \in \prod_{i=1}^{a}\Lambda^{e_1}}f_\mu(D_v)$. By Theorem~\ref{general Cuntz-Krieger single vertex k-graph}, in order to prove that $\pi$ is injective we only need to show that the spectrum of $\pi(W)$ in $C^*(\pi(W))$ contains the unit circle. By \cite[Proposition~3.11]{Katsura:CJM03}, it suffices to construct an expectation $\Phi:C^*(\pi(W)) \to \mathbb{C} \cdot \pi(1_{C^*(\Lambda)})$ such that $\Phi(\pi(1_{C^*(\Lambda)}))=\pi(1_{C^*(\Lambda)}),\Phi(\pi(W^n))=0$ for all $n \in \mathbb{Z} \setminus \{0\}$. Fix $\{z_n\}_{n \in \mathbb{Z}} \subset \mathbb{C}$ with at most finitely many nonzero. Let $\mathcal{F}:=\{0 \neq n \in \mathbb{Z}:z_n \neq 0\}$. By the assumption of the theorem, there exist $\mu_0 \in \prod_{i=1}^{a}\Lambda^{e_1}$ and a measurable subset $E$ of $f_{\mu_0}(D_v)$ such that $\eta(E) \neq 0$ and $(f_{\mu_0} \circ f_{h(\mu_0)}^{-1})^{n}(E) \cap E\stackrel{\eta-a.e.}{=}\emptyset$ for all $n \in \mathcal{F}$.  Take an arbitrary function $\phi \in L^2(X,\eta)$ with $\Vert\phi\Vert = 1$ and $\supp(\phi) \stackrel{\eta-a.e.}{\subset}E$. Then
\begin{align*}
\Big\Vert \sum_{n \in \mathbb{Z}}z_n\pi(W^n) \Big\Vert^2 &=\Big\Vert z_0 \pi(1_{C^*(\Lambda)})+\sum_{n \in \mathcal{F}}z_n\pi(W^n)\Big\Vert
\\&\geq \Big\Vert z_0\phi \chi_{D_v} + \sum_{n \in \mathcal{F}}z_n\pi(W^n)(\phi) \Big\Vert^2
\\&=\int_X  \Big\vert z_0\phi\chi_{D_v} + \sum_{n \in \mathcal{F}}z_n\pi(W^n)(\phi) \Big\vert^2 \, \mathrm{d}\eta
\\&\geq\int_E  \Big\vert z_0\phi + \sum_{n \in \mathcal{F}}z_n\pi(W^n)(\phi) \Big\vert^2 \, \mathrm{d}\eta
\\&=\int_E  \Big\vert z_0\phi + \sum_{n \in \mathcal{F}}z_n\pi(s_{h(\mu_0)}s_{\mu_0}^*)(\phi) \Big\vert^2 \, \mathrm{d}\eta
\\& = \int_E\Big\vert z_0 \phi \Big\vert^2\,\mathrm{d}\mu
\\&=\vert z_0\vert^2.
\end{align*}
So we get the required expectation $\Phi$ and hence $\pi$ is injective. \fim

In the following we modify the construction of the branching systems in Theorem~\ref{existenceofabranchingsystem} and we obtain a branching system for each periodic single-vertex $2$-graph so that the associated representation is faithful.

\begin{exemplo}\label{construction of branching system for k-graph}
Let $\Lambda$ be a periodic single-vertex $2$-graph with $\vert \Lambda^{e_1}\vert, \vert \Lambda^{e_2}\vert\geq 2$. Let $X:=[0,1] \times \Lambda^{\infty}$. Define $D_v:=X$. For $e \in \Lambda^{e_1}$, define $R_e:=[0,1] \times e\Lambda^{\infty}$, define $F_e:D_v \to R_e$ by $F_e(t,x):=(t^2,ex)$ (see Lemma~\ref{sigma is a local homeo}). For $f \in \Lambda^{e_2}$, define $R_f:=[0,1] \times f\Lambda^{\infty}$, define $F_f:D_v \to R_f$ by $F_f(t,x):=(\sqrt{t},fx)$. Then $\{D_v,R_\mu,F_\mu\}_{\mu \in \bigcup_{i=1}^{2}\Lambda^{e_i}}$ is a $\Lambda$-branching system. By Theorem~\ref{a criterion of faithful rep}, the induced representation is faithful.
\end{exemplo}

We finish this section by building a branching system on $\R^2$ for a periodic single-vertex $2$-graph such that the associated representation is faithful.

\begin{exemplo}\label{exemploinjetor}

Consider the flip C*-algebra from the 2-colored graph of \cite[Example 4.3]{DavidsonYang},

\centerline{
\setlength{\unitlength}{2cm}
\begin{picture}(1,0.7)
\color{blue}
\put(0.5,0){\qbezier(-0.05,0)(-1,-1)(-1,0)}
\put(0.5,0){\qbezier(-0.05,0)(-1,1)(-1,0)}
\put(-0.3,0.44){$>$}
\put(-0.75,0){$f_2$}
\put(0.5,0){\qbezier(-0.05,0)(-0.8,-0.4)(-0.8,0)}
\put(0.5,0){\qbezier(-0.05,0)(-0.8,0.4)(-0.8,0)}
\put(-0.2,0.142){$>$}
\put(-0.47,0){$f_1$}
\color{red}
\put(0.5,0){\qbezier(0.05,0)(0.8,0.4)(0.8,0)}
\put(0.5,0){\qbezier(0.05,0)(0.8,-0.4)(0.8,0)}
\put(0.5,0){\qbezier(0.05,0)(1,1)(1,0)}
\put(0.5,0){\qbezier(0.05,0)(1,-1)(1,0)}
\put(1.2,0.445){$<$}
\put(1.55,0){$e_2$}
\put(1,0.14){$<$}
\put(1.34,0){$e_1$}
\color{black}
\put(0.5,0){\circle*{0.08}}
\put(0.45,0.1){$v$}
\end{picture}}
\vspace{1.5 cm}

with the factorization rule $e_if_j=f_ie_j$, for $i,j\in \{1,2\}$. For this 2-graph, here called $\Lambda$, we obtain a $\Lambda$-branching system on the measure space $\left([0,1]\times [-1,1], \eta\right)$ in the sense of Theorem \ref{a criterion of faithful rep}, where $\eta$ is the Lebesgue measure in $\R^2$.
\end{exemplo}

\begin{proposicao}
Let $D_v=[0,1]\times[-1,1]$, $R_{e_1}=R_{f_1}=[0,1]\times [0,1]$, and $R_{e_2}=R_{f_2}=[0,1]\times [-1,0]$. Define the maps $f_{e_1}, f_{f_1}:D_v\rightarrow R_{e_1}$ by $f_{e_1}((x,y))=(x^2,\frac{y}{2}+\frac{1}{2})$ and $f_{f_1}((x,y))=(\sqrt{x},\frac{y}{2}+\frac{1}{2})$, and the maps $f_{e_2},f_{f_2}:D_v\rightarrow R_{f_2}$ by $f_{e_2}((x,y))=(x^2,\frac{y}{2}-\frac{1}{2})$ and $f_{f_2}((x,y))=(\sqrt{x},\frac{y}{2}-\frac{1}{2})$. Then the representation of $C^*(\Lambda)$ arising from this branching system is faithful.
\end{proposicao}

\demo
It is easy to see that conditions all the conditions of Definition \ref{branchsystem fin aligned case} are satisfied.
Let $\pi:C^*(\Lambda)\rightarrow B(L^2([0,1]\times [-1,1],\eta))$be the *-homomorphism induced by this $\Lambda$-branching system. Note that each cycline pair is of the form $(e_i, f_i)$. So, to apply Theorem \ref{a criterion of faithful rep}, it is enough to show that, for each finite set $\mathcal{F}\subseteq \Z\setminus\{0\}$, there exists a subset $E\subseteq f_{e_1}(D_v)$ with $\eta(E)\neq 0$ and $(f_{e_1}\circ f_{f_1}^{-1})^k(E)\cap E\stackrel{\eta-a.e.}{=}\emptyset$  for each $k\in\mathcal{F} $. Note that $(f_{e_1}\circ f_{f_1}^{-1})(x,y)=(x^4,y)$.

Let $E=[\frac{1}{4}, \frac{1}{2}]\times [0,1]$. Then $(f_{e_1} \circ f_{f_1}^{-1})^k(E)\cap E\stackrel{\eta-a.e.}{=}\emptyset$ for each $k\in \Z\setminus\{0\}$, and hence, by Theorem~\ref{a criterion of faithful rep}, $\pi$ is faithful.
\fim

\section*{Acknowledgments}


The authors would like to thank Dr. Jonathan Brown and  Prof. Dilian Yang for valuable discussions regarding the theory of $k$-graph C*-algebras. The second author also thanks the rest of the authors for the continuation of the collaboration.

\vspace{1.5pc}
\begin{center}
Daniel Gon\c{c}alves (daemig@gmail.com) and Danilo Royer (danilo.royer@ufsc.br)\\

Departamento de Matem\'{a}tica - Universidade Federal de Santa Catarina, Florian\'{o}polis, 88040-900, Brazil

\vspace{1.5pc}

Hui Li (lihui8605@hotmail.com)\\

Research Center for Operator Algebras and Shanghai Key Laboratory of Pure Mathematics and Mathematical Practice, Department of Mathematics, East China Normal University, 3663 Zhongshan North Road, Putuo District, Shanghai 200062, China \\ Department of Mathematics $\&$ Statistics, University of Windsor, Windsor, ON N9B 3P4, CANADA

\end{center}

\end{document}